\documentclass{article}
\usepackage{amsmath}   
\usepackage{amssymb}   
\usepackage{accents}
\usepackage{graphicx}
\usepackage{color}
\usepackage{tikz}
\usepackage{tikz-cd}
\usetikzlibrary{shapes.geometric,fit,trees}
\usepackage{forest}
\usepackage{paralist}
\usepackage{enumerate}    
\usepackage{amsthm}  
\usepackage[english]{babel}  
\usepackage{textgreek}
\usepackage[normalem]{ulem}
\usepackage{appendix}
\usepackage[shortlabels]{enumitem}

\theoremstyle{theorem}            
\newtheorem{theorem}{Theorem}
\theoremstyle{definition}           

\newtheorem{definition}[theorem]{Definition}

\newtheorem{proposition}[theorem]{Proposition}
\newtheorem{lemma}[theorem]{Lemma}
\newtheorem{corollary}[theorem]{Corollary}

\newtheorem{question}{Question}

\newtheorem{notation}[theorem]{Notation}

\newtheorem*{remark}{Remark}

\newcommand{\ignore}[1]{}
\newcommand{\term}{\textbf} 
\newcommand{\concat}{^\smallfrown} 
\renewcommand{\models}{\vDash} 
\newcommand{\notmodels}{\nvDash} 
\newcommand{\forces}{\Vdash} 
\newcommand{\0}{\varnothing} 
\renewcommand{\bar}{\overline} 
\renewcommand{\phi}{\varphi} 
\renewcommand{\epsilon}{\varepsilon} 

\DeclareMathOperator{\range}{range} 

\def\N{\mathbb{N}} 
\def\R{\mathbb{R}} 
\def\Cantor{2^\N} 
\def\Baire{\N^\N} 
\def\Ord{\textbf{Ord}} 

\def\KP{\mathsf{KP}}
\def\ZF{\mathsf{ZF}}
\def\ZFC{\mathsf{ZFC}}

\newcommand{\lo}{\textup{LO}}
\newcommand{\loiso}{E_{\textup{LO}}}

\def\acts{\curvearrowright}

\DeclareMathOperator{\Aut}{Aut}
\def\Fw{F_{\omega_1}}
\def\Fwrel{F_{\omega_1}^\text{rel}}
\def\O{\mathcal{O}}

\DeclareMathOperator{\proj}{proj}
\def\Mod{\mathsf{Mod}}
\def\lang{\mathcal{L}}

\title{Borel graphable equivalence relations}
\author{Tyler Arant\footnote{University of California, Los Angeles}, Alexander S. Kechris\footnote{California Institute of Technology; Research partially supported by NSF Grant DMS-1950475} and Patrick Lutz\footnote{University of Michigan, Ann Arbor; Research partially supported by NSF Grant DMS-203072}}
\date{\today}

\begin{document}

\maketitle

\vspace{-10pt}

\begin{abstract}
This paper is devoted to the study of analytic equivalence relations which are \textit{Borel graphable}, i.e. which can be realized as the connectedness relation of a Borel graph. Our main focus is the question of which analytic equivalence relations are Borel graphable. First, we study an equivalence relation arising from the theory of countable admissible ordinals and show that it is Borel graphable if and only if there is a non-constructible real. As a corollary of the proof, we construct an analytic equivalence relation which is (provably in $\ZFC$) not Borel graphable and an effectively analytic equivalence relation which is Borel graphable but not effectively Borel graphable. Next, we study analytic equivalence relations given by the isomorphism relation for some class of countable structures. We show that all such equivalence relations are Borel graphable, which implies that for every Borel action of $S_\infty$, the associated orbit equivalence relation is Borel graphable. This leads us to study the class of Polish groups whose Borel actions always give rise to Borel graphable orbit equivalence relations; we refer to such groups as \textit{graphic groups.} We show that besides $S_\infty$, the class of graphic groups includes all connected Polish groups and is closed under countable products. We finish by studying structural properties of the class of Borel graphable analytic equivalence relations and by considering two variations on Borel graphability: a generalization with hypergraphs instead of graphs and an analogue of Borel graphability in the setting of computably enumerable equivalence relations.
\end{abstract}

\bigskip

\noindent{\bf (A)} Every undirected graph $G$ on a set $X$ induces a connectedness equivalence relation, $E_G$, defined by
\[
xE_Gy \iff \text{there is a path from $x$ to $y$ in $G$} \qquad
(x, y\in X).
\]
If $X$ is a Polish space and the graph $G$ is Borel, then the connectedness equivalence relation $E_G$ is analytic. 
We will call analytic equivalence relations which arise in this way Borel graphable (a notion first studied in the PhD thesis of the first author~\cite{arant2019effective}).

\begin{definition}
An equivalence relation $E$ on a Polish space $X$
is \term{Borel graphable} if $E$ is the connectedness equivalence
relation of a Borel graph $G$ on $X$.  The graph $G$ is called
a \term{Borel graphing} of $E$.
\end{definition}

Note that we can similarly define the notion of a $\Phi$-graphable equivalence relation for any pointclass $\Phi$.  For instance, in the effective theory, we can consider the class of $\Delta^1_1$-graphable equivalence relations.

The main goal of this paper is to understand which analytic equivalence relations are Borel graphable. Thus we will begin by briefly reviewing the theory of analytic equivalence relations.

\medskip\noindent
{\bf (B)} There are a great number of important examples of analytic equivalence relations occurring in several areas of mathematics, often associated with some notion of isomorphism. We mention a few examples that we will study in this paper.

\begin{enumerate}[(i)]
\item The isomorphism relation on the countable structures with universe $\N$ of a fixed countable first-order language. Similarly, the biembeddability relation on such structures.

\item The equivalence relation $E^\Gamma_X$ induced by a Borel action of a Polish group $\Gamma$ on a Polish space $X$.
\end{enumerate}

By Burgess' Theorem (see \cite{burgess1978equivalences}) every analytic equivalence relation has countably many, $\aleph_1$ many, or continuum many equivalence classes. The next examples are canonical analytic equivalence relations with exactly $\aleph_1$ many classes.

\begin{enumerate}[resume*]
\item On $2^\N$, the equivalence relation
\[
x F_{\omega_1} y \iff \omega_1^x = \omega_1^y,
\]
where $\omega_1^x$ is the Church-Kleene ordinal of $x$.

\item On the space LO of linear orderings on $\N$, the relation
\[
\leq_1 E_{\omega_1} \leq_2 \iff (\leq_1, \leq_2 \ \in {\rm WO}\  \&  \ \leq_1 \cong \leq_2) \vee (\leq_1, \leq_2  \ \notin {\rm WO}),
\]
where WO is the set of wellorderings on $\N$.

\item (see \cite{lauchli1988linearly}) On the space LO, the relation
\[
\leq_1 \cong_p \leq_2 \iff \leq_1 \textrm{is piecewise isomorphic to} \leq_2,
\]
where, for a relational language, two structures $M,N$ in this language are {\bf piecewise isomorphic} if there is a bijection $f\colon M \to N$ and a partition $A_0, A_1, \dots, A_{n-1}$ of $M$, such that for all $i$, $f|A_i$ is an isomorphism of $M|A_i$ with $N|f(A_i).$
\end{enumerate}

One can organize analytic equivalence relations under the hierarchical order of Borel reducibility. Given such relations $E,F$ on $X,Y$, resp., a Borel reduction of $E$ to $F$ is a Borel map $f\colon X\to Y$ such that 
\[
x_1 E x_2 \iff f(x_1) F f(x_2).
\]
If such $f$ exists, we say that $E$ is {\bf Borel reducible} to $F$, in symbols $E\leq_B F$. An analytic equivalence relation is called universal if all other analytic equivalence relations can be Borel reduced to it; it is known that such universal equivalence relations exist.

\medskip\noindent
{\bf (C)} We will now give a brief overview of the main results of this paper. In Section \ref{sec:basics}, we establish the basic properties of Borel graphability and prove the Borel graphability of some motivating examples.  We also use these examples to highlight some general notions and methods that will
be used throughout. In particular, we show that there is a universal analytic equivalence relation that is Borel graphable but that not all analytic equivalence relations are Borel graphable. We also introduce the notion of the {\bf Borel witness coding property} for an analytic equivalence relation and show that relations with this property have Borel graphings of diameter 2. Using this we show that isomorphism of linear orderings of $\N$, as well as biembeddability and piecewise isomorphism of linear orderings of $\N$, and $E_{\omega_1}$ all admit Borel graphings of diameter 2.

In Section \ref{sec:admiss}, we study the analytic equivalence relation $F_{\omega_1}$. We show the surprising result that the Borel graphability of this equivalence relation is equivalent to the existence of a non-constructible real.  By taking a direct sum of the relativized versions of this equivalence relation, we encounter a complicated example of an analytic equivalence relation which is not Borel graphable.

In Section \ref{sec:isomorphism}, we establish that for any countable language, the isomorphism relation on the space of countable models is Borel graphable (extending the results from Section~\ref{sec:basics} on linear orders).  We also show that for any countable language, the associated piecewise isomorphism and biembeddability equivalence relations are Borel graphable. 
  
In Section \ref{sec:graphic}, we establish several results about Borel graphability of equivalence relations arising from Borel actions of Polish groups. In particular, we consider the class of \term{graphic} Polish groups---Polish groups for which every Borel action of the group on a Polish space gives rise to a Borel graphable orbit equivalence relation. We show that all connected Polish groups are graphic and that all automorphism groups of countable ultrahomogeneous structures in a finite relational language (including $S_\infty$) are graphic. We also prove a few closure properties for the class of graphic Polish groups.

In Section \ref{sec:closure_properties} we examine structural properties of the class of Borel graphable equivalence relations, including closure properties. Some graphability results in previous sections make use of coding properties for equivalence relations that are weakenings of the Borel witness coding property from Section~\ref{sec:basics}. In Section \ref{sec:borelcoding} we explore the relationships between these coding properties and Borel graphability.

The final sections are about some variations of the notion of Borel graphability. Borel graphability can be generalized by using hypergraphs, and we explore this notion in Section \ref{sec:hypergraph}. In Section \ref{sec:ceers}, we explore the analogue of Borel graphability in the setting of computably enumerable equivalence relations and establish that any computably enumerable equivalence relation all of whose classes are infinite is computably graphable (a related result was independently obtained by Uri Andrews and coauthors). 

There is also an Appendix, which consists of an exposition of Kumabe-Slaman forcing, a key tool used in Section \ref{sec:admiss}.

\subsection*{Acknowledgements}

Thanks to Gabe Goldberg for pointing out the absoluteness argument in Section~\ref{sec:admiss}, to Jan Greb\'ik for several helpful conversations and to an anonymous reviewer for pointing out a simplified proof of Theorem~\ref{thm:noneffective}.  A special thanks to Andrew Marks for many useful
comments and insights when the notion of Borel graphability was first considered.


\section{Basic properties and examples}
\label{sec:basics}

As mentioned above, in this section we will establish some basic properties of Borel graphability and see a few examples that demonstrate some of the techniques we will use throughout the paper. We begin with the following proposition which establishes that the question of Borel graphability is only interesting for analytic, non-Borel equivalence relations which are not countable (i.e.\ which have at least one uncountable class).

\begin{proposition}
\label{trivialobsv}
Let $E$ be an equivalence relation on a Polish space $X$.
\begin{enumerate}
    \item If $E$ is Borel, then $E$ is Borel graphable.
    \item If $E$ is Borel graphable, then $E$ is analytic.
    \item If $E$ is countable and Borel graphable, then $E$ is Borel.
\end{enumerate}
\end{proposition}

\begin{proof}
(1) and (2) are trivial.  We prove (3). Note that if $X$ is countable then $E$ is trivially Borel, so we may assume $X$ is uncountable. Since any two uncountable Polish spaces are Borel isomorphic (see Theorem 15.6 of~\cite{kechris1995classical}), we may assume without loss of generality that $X = \Cantor$.  Suppose $E$ is countable and let $G$ be a Borel graphing of $E$.  Let $P \subseteq X^3$ be the set defined by
\[
    (x, y, z) \in P \iff z \text{ codes a path\footnotemark~in $G$ from $x$ to $y$}.
\]
\footnotetext{Here is more precisely what we mean by ``$z$ codes a path in $G$.'' Fix a Borel bijection $\phi \colon \Cantor \to (\Cantor)^{< \N}$ (which exists since all uncountable Polish spaces are Borel isomorphic). By ``$z$ codes a path in $G$,'' we mean that if $\phi(z) = (z_0, z_1, \ldots, z_n)$ then $z_0, z_1, \ldots, z_n$ forms a path in $G$.}
Since $E$ is countable, it follows that $G$ is locally countable and hence that $P$ has countable sections---i.e.\ for any $x, y \in X$, the set $P_{x, y} = \{z \mid (x, y, z) \in P\}$ is countable. Note that
\[
    x E y \iff (\exists z)~(x, y, z) \in P.
\]
It follows from the Lusin-Novikov Theorem (see Theorem 18.10 of~\cite{kechris1995classical}) that $E$ is Borel.
\end{proof}

It is natural to ask whether every analytic equivalence relation is Borel graphable. Using item (3) of Proposition \ref{trivialobsv} above, we can easily construct examples showing that this is not the case.

\begin{proposition}
\label{prop:firstnonexample}
There is an analytic equivalence relation which is not Borel graphable.
\end{proposition}

\begin{proof}
Let $A$ be an analytic, non-Borel
subset of a Polish space $X$ and define an equivalence relation 
$E$ on $X\times \{0, 1\}$ by
\[
(x, i) E(y, j) \iff ((x, i)=(y, j)  \text{ or } x=y\in A).
\]
In other words, for every $x \in A$, $(x, 0)$ and $(x, 1)$ form a single equivalence class of size $2$ and for every $x \notin A$, $(x, 0)$ and $(x, 1)$ form two separate equivalence classes, each of size $1$.  Since all equivalence classes of $E$ are countable, Proposition~\ref{trivialobsv} implies that $E$ is Borel graphable if and only if $E$ is Borel. But since $A$ is not Borel, neither is $E$.
\end{proof}

On the other hand, there are complicated equivalence relations which are Borel graphable, including a universal one.

\begin{proposition}
    There is a universal analytic equivalence relation
    which is Borel graphable.  
\end{proposition}

\begin{proof}
    Let $E$ be a universal analytic equivalence relation on $\Baire$. 
    Pick a closed set $C\subseteq (\Baire)^3$ such that
    $E$ is the projection of $C$.  Define a graph $G$ on the space $\Baire \cup C$ so that $x\in \Baire$ is
    only adjacent to elements $(y, y', z)$ in $C$ with $x\in\{y, y'\}$.
    Then, it is easy to check that the connectedness relation of
    $G$ is the desired universal equivalence relation.
\end{proof}

We now turn to the question of what sort of operations preserve Borel graphability. As a consequence of the two propositions above, we have:

\begin{corollary}
\label{cor:borelreducibility}
Borel graphability is not closed downwards under Borel reducibility.
\end{corollary}

Although Borel reducibility does not preserve Borel graphability, there is a stronger sort of reducibility that does. Recall that a reduction $f$ from $E$ to $F$ is \term{invariant} if the image of $f$ is $F$-invariant, i.e., a union of $F$-classes.

\begin{proposition}
\label{prop:invariantreducability}
Borel graphability is closed downward under invariant Borel reducibility.
\end{proposition}

\begin{proof}
The invariance of the reduction allows one to
transfer over the graphing in a way that preserves connectedness.
\end{proof}

As an aside, we now mention an interesting basic question about Borel graphability. Recall that the \term{diameter} of a graph is the maximum distance between any two connected vertices.  In many of the results throughout this paper, we will explicitly state the diameter of the Borel graphing obtained in the proof.  Very often, we obtain a graphing of diameter $2$, but we do not currently know whether this is always possible---that is, whether every Borel graphable equivalence relation has a Borel graphing with diameter $2$.

\begin{question}
    Is there a Borel graphable equivalence relation $E$
    such that every Borel graphing of $E$ has diameter
    $\geq 3$?
\end{question}

\begin{remark}
After this paper was written, the third author answered this question by proving that there is an analytic equivalence relation $E$ such that $E$ has a Borel graphing with diameter $4$, but no graphing with diameter $2$ (in fact, $E$ is the equivalence relation defined in the proof of Theorem~\ref{thm:separation} below). However, the proof does not seem to resolve most variations on the above question. For example, it does not seem to help determine whether there is a Borel graphable equivalence relation which has no Borel graphing of finite diameter.
\end{remark}

Next, we observe that there are no local obstructions
to Borel graphability, in the sense that no
single class on its own can prevent a Borel graphing.

\begin{theorem}
\label{analyticset}
Let $X$ be a Polish space and let $A\subseteq X$ be an analytic set.
Then there is a Borel graph $G\subseteq A^2$ on $X$ of diameter $2$ such that $A$
is a connected component of $G$. 
\end{theorem}

\begin{proof}
By the perfect set theorem for analytic sets, either $A$ is countable or $A$ contains a perfect set. We will handle each of these two cases separately.

\medskip\noindent\textit{Case 1: $A$ is countable.} This means $A$ is Borel, so we can take $G = \{(x,y) \in A^2\colon x\not= y\}$.

\medskip\noindent\textit{Case 2: $A$ is uncountable.} 
By the perfect set theorem for analytic sets (see Exercise 14.13 of~\cite{kechris1995classical} or Theorem 2C.2 of~\cite{moschovakis2009}), $A$ contains a Cantor set, i.e.\ a subset which is homeomorphic to $\Cantor$. Let $P \subseteq A$ be such a set and let $\phi \colon P \to \Cantor$ be a homeomorphism. It is a standard fact that since $A$ is analytic, it is the projection of a $G_\delta$ subset of $X\times \Cantor$ (see Exercise 14.3 of~\cite{kechris1995classical}). Let $B \subseteq X\times \Cantor$ be such a set, i.e.
\[
    x\in A \iff \exists y\in \Cantor\, (x,y)\in B.
\]
Now define $G'$ as follows
\[
    (x,y)\in G' \iff [y \in P \  \&  \ (x, (\phi(y))_0) \in B] \vee [y\in P  \ \&  \ (x, (\phi(y))_1) \in B],
\]
where $z\mapsto ((z)_0, (z)_1)$ is a homeomorphism of $\Cantor$ with $(\Cantor)^2$.

Clearly $G'$ is $G_\delta$. Let $G$ be the symmetrization of $G'$ minus the diagonal, so that $G$ is an undirected, simple graph which is also $G_\delta$ and $G\subseteq A^2$.

Let now $x, z\in A, x\not= z$. We will find $y$ such that $(x,y)\in G', (z,y)\in G'$, which completes the proof. Indeed, let $a,b$ be such that $(x,a) \in B, (z,b)\in B$ and let $y\in P$ be such that $(\phi(y))_0 = a, (\phi(y))_1 = b$. Then clearly $(x,y)\in G', (z,y)\in G'$.
\end{proof}

We remark that the proof of the previous theorem contains
the core idea of many of the arguments for Borel graphability.
Often, a graphing $G$ for an equivalence relation $E$ can be
defined by placing an edge between equivalent $x$ and $y$
when $y$ codes enough information, e.g., $y$ codes a witness
to $xEy$.  The following definitions provide a useful way
to formalize this idea.  

\begin{definition}
\label{def:unfolding}
Given an analytic equivalence relation $E$ on $X$, an {\bf unfolding} of $E$ is a Borel set $C\subseteq X^2\times \Cantor$ such that
\begin{enumerate}
    \item $\forall z,~ (x,x,z) \in C$
    \item $(x,y,z) \in C \iff (y,x,z)\in C$
    \item $x E y \iff \exists z (x,y,z) \in C$.
\end{enumerate}
We call any $z$ such that $(x,y,z) \in C$ a {\bf witness} to $x E y$.\footnote{Note that the definition of unfolding implies that
any witness to $xEy$ is a witness to $yEx$ and that every
$z\in\Cantor$ is a witness to $xEx$.}
\end{definition}

\begin{definition}
We say that an analytic equivalence relation $E$ on $X$ has the \term{Borel witness coding property} if there is an unfolding $C$ of $E$ and a Borel function $f:X\rightarrow\Cantor$ such that for all $x \in X$ and $a \in \Cantor$, there is some $y$ which is $E$-equivalent to $x$ such that $f(y)$ computes both $a$ and a witness to the equivalence $x E y$.
\end{definition}

Typically, the function $f(y)$ in the definition
will just be some small number of jumps of (a code for) $y$, perhaps relative to some parameter. In some cases, $f(y)$ will even just be a code for $y$ itself, with no jumps or parameters required.\footnote{We have already used the words ``code'' and ``coding'' a few times without being precise about what we mean. There are many reasonable ways to construct codes for points in a Polish space, but we will now specify one such way for concreteness. Every Polish space $X$ is homeomorphic to a $G_\delta$ subset of the Hilbert cube, $[0, 1]^\N$ (see Theorem 4.14 of~\cite{kechris1995classical}) Furthermore, we can naturally associate each element $(a_i)_i$ in $[0, 1]^\N$ to an element of $\Cantor$ by, say, taking the join of the characteristic functions of the Dedekind cuts of the $a_i$. By a \textit{code} for $y \in X$, we mean an element of $\Cantor$ which is naturally associated to an element of $[0, 1]^\N$ which maps to $y$ via some fixed homeomorphism between a $G_\delta$ subset of $[0, 1]^\N$ and $X$. The general idea is that a code for $y$ is an element of $\Cantor$ which ``knows'' about which basic neighborhoods $y$ belongs to.}

\begin{proposition}
    If an analytic equivalence relation $E$ 
    has the Borel witness coding property, then $E$
    is Borel graphable by a graph of diameter 2.
\end{proposition}

\begin{proof}
    We can assume that $E$ is on $X=\Cantor$. Let $f:X\rightarrow \Cantor$ be a function
    witnessing that the Borel witness coding property holds for $E$ with respect to the unfolding $C$.
    Define first $G'\subseteq X^2$ by
    \begin{multline*}
    (x, y)\in G' \iff f(x) \ \text{computes some $z\in X$}\\
    \text{and computes witnesses that $xEz$ and $zEy$.}
    \end{multline*}
    Next, let $G$ be the graph defined by symmetrizing $G'$ and removing
    the diagonal.  Clearly, $G$ is Borel, and we show
    that it is a graphing of $E$.\footnote{Note that quantifying
    over the reals computable from $f(x)$ is quantification over
    a countable set, hence it is a Borel operation.  This fact, and
    an analogous version for hyperarithmetical reduction, will be
    used many times in the sequel.}  

    Suppose $xEy$.  Use the witness coding property
    to obtain a $z \in X$ such that $f(z)$ computes $x$ and computes
    witnesses to $xEy$ and $xEz$.  Immediately, it is clear
    that $(z, y)\in G$.  Also, $(x, z)\in G$
    since $f(z)$ computes $x$ and witnesses to $xEx$ and $xEz$.
\end{proof}

One useful feature of the Borel witness coding property is that it is often easy to transfer from an equivalence relation $E$ to a related equivalence relation $E'$; thus, we can immediately infer Borel graphability of $E'$, even in situations where Borel graphability of $E'$ cannot be obviously inferred from Borel graphability of $E$. The propositions below give an example of this. 

\begin{proposition}
\label{subsetwitness}

  Let $E$ and $E'$ be equivalence relations on a Polish space
    $X$.  If $E\subseteq E'$ and $E$ has the Borel witness
    coding property, then $E'$ has the Borel witness coding property.  In particular, $E'$ is Borel
    graphable by a graph with diameter 2.
\end{proposition}

\begin{proof}
    Pick unfoldings $C,C'$ for $E,E'$, resp., and a Borel function $f$ which verifies the Borel witness coding property for $E$ with respect to $C$. Since $E\subseteq E'$, $C\cup C'$ is an unfolding of $E'$ and $f$, together with $C\cup C'$, demonstrates the Borel witness coding property for $E'$. 
\end{proof}

We will now see a simple example of how to establish the
Borel witness coding property, and hence Borel graphability.
Let $\lo$ denote the Polish space of linear orders with domain $\N$,
and let $\loiso$ be the isomorphism equivalence relation
on $\lo$ (for a more detailed definition of $\lo$ and a discussion of its properties, see Sections 16C and 27C of~\cite{kechris1995classical}).

\begin{theorem}
\label{thm:loiso}
    $\loiso$ has the Borel witness coding property.  
    In particular, $\loiso$ is Borel graphable with
    a graph with diameter $2$.
\end{theorem}

\begin{proof}
In order to prove that $\loiso$ has the Borel witness coding property, we first need to pick an unfolding of $\loiso$. In later proofs, we will typically leave this step implicit as it will usually be fairly clear what unfolding is intended, but here we describe it explicitly to illustrate our definitions. First recall that $\Cantor$ is homeomorphic to $2^{\N\times \N}$ and thus we may view an element of $\Cantor$ as a subset of $\N\times \N$. Given that, our unfolding $C \subseteq \lo\times \lo \times \Cantor$ will consist of all tuples $(\leq_1, \leq_2, x)$ such that either $\leq_1 = \leq_2$ or $x$, thought of as a subset of $\N\times \N$, is the graph of an isomorphism from $\leq_1$ to $\leq_2$ or from $\leq_2$ to $\leq_1$ (recall that $\leq_1$ and $\leq_2$ are both linear orders on $\N$). 

Now let $\leq_1\in\lo$ and $a\in \Cantor$ be given. To show that $\loiso$ has the Borel witness coding property, it is enough to build a linear order $\leq_2$ on $\N$ which is isomorphic to $\leq_1$ and which computes both $a$ and an isomorphism between $\leq_1$ and $\leq_2$.

The main idea of the construction is that information can be coded into $\leq_2$ by choosing, for each $n$, whether to make $2n <_2 (2n + 1)$ or $(2n + 1) <_2 2n$. Using this idea, we can build $\leq_2$ which is both isomorphic to $\leq_1$ but also codes a lot of information by more-or-less copying $\leq_1$, but sometimes swapping adjacent elements (where ``adjacent'' is in the sense of the standard ordering on $\N$, not in the sense of $\leq_1$). 

Also, instead of constructing $\leq_2$ directly, we will first construct a bijection $f\colon \N \to \N$ and then define $\leq_2$ so that $f$ is an isomorphism from $\leq_2$ to $\leq_1$---in other words, we will set $n \leq_2 m$ if and only if $f(n) \leq_1 f(m)$. For each $n$, $f$ will either map $2n$ and $2n + 1$ to themselves or it will swap them, depending on the value of the $n^\text{th}$ bit of information to be encoded.  

We will now explain this in more detail. Fix a code $b \in \Cantor$ for $\leq_1$ and define $c = a\oplus b$; we will encode $c$ into $\leq_2$ (as we will see, this is enough to give us everything we need). For each $n$, say that $c$ and $\leq_1$ \term{agree at $n$} if either $c(n) = 0$ and $2n <_1 2n + 1$ or $c(n) = 1$ and $2n + 1 <_1 2n$ (i.e.,\ if $\leq_1$ already correctly codes the $n^\text{th}$ bit of $c$). Now define $f \colon \N \to \N$ by
\begin{align*}
f(2n) =
\begin{cases}
2n &\text{if $c$ and $\leq_1$ agree at $n$.}\\
2n + 1 &\text{otherwise.}
\end{cases}\\
f(2n + 1) =
\begin{cases}
2n + 1 &\text{if $c$ and $\leq_1$ agree at $n$.}\\
2n &\text{otherwise.}
\end{cases}
\end{align*}
In other words, if $c$ and $\leq_1$ agree at $n$ then $f$ maps $2n$ and $2n + 1$ to themselves and otherwise it swaps them. Finally, define $\leq_2$ by setting $n \leq_2 m$ if and only if $f(n) \leq_1 f(m)$.

We now check that $\leq_2$ has the desired properties. First note that $\leq_2$ is isomorphic to $\leq_1$, as witnessed by $f$. Next, note that $c$ can be computed from $\leq_2$: to check if the $n^\text{th}$ bit of $c$ is $0$, we can just check if $2n <_2 2n + 1$. Thus $\leq_2$ computes both $a$ and $\leq_1$. Given that $\leq_2$ computes $\leq_1$, it is easy to see that it computes $f$ as well: to check whether $f$ keeps $2n$ and $2n + 1$ the same or swaps them, simply check if $\leq_1$ and $\leq_2$ agree on the order of $2n$ and $2n + 1$ or not. Thus $\leq_2$ computes both $a$ and an isomorphism between $\leq_2$ and $\leq_1$, as required.
\end{proof} 

\begin{corollary}
\label{cor:lo}    
    The following equivalence relations on $\lo$ are all
    Borel graphable by graphs of diameter $2$.
    \begin{enumerate}
    \item Biembeddability of linear orders.
    \item Piecewise isomorphism of linear orders.
    \item $E_{\omega_1}$.
    \end{enumerate}
\end{corollary}

\begin{proof}

    Follows immediately from Proposition \ref{subsetwitness}
    and Theorem \ref{thm:loiso}.
\end{proof}

\begin{remark}
Theorem~\ref{thm:loiso} and Corollary~\ref{cor:lo} will be generalized quite a bit in Section~\ref{sec:isomorphism} and even further in Section~\ref{sec:graphic}.
\end{remark}


\section{Countable admissible ordinals}
\label{sec:admiss}

In this section, we will focus on the equivalence relation $\Fw$, which arises from computability theory. In order to make this section more accessible, we will begin by reviewing some of the definitions involved. Let $\omega_1^{CK}$ denote the least ordinal with no computable presentation (sometimes referred to as the \term{Church-Kleene ordinal}) and for any real $x \in \Cantor$, let $\omega_1^x$ denote the least ordinal which has no presentation computable from $x$. Note that for every $x$, $\omega_1^x$ is countable and $\omega_1^{0} = \omega_1^{CK}$.

Now let $\Fw$ be the equivalence relation on $\Cantor$ defined by
\[
    x \Fw y \iff \omega_1^x = \omega_1^y.
\]
$F_{\omega_1}$ is sometimes referred to as the \term{countable admissible ordinal equivalence relation} because ordinals of the form $\omega_1^x$ are exactly the countable admissible ordinals.\footnote{An admissible ordinal is an ordinal $\alpha$ such that $L_\alpha \models \KP$, where $\KP$ denotes Kripke-Platek set theory, a weak fragment of $\ZF$ studied in higher computability theory. We will not make any use of this characterization.} It is known that $F_{\omega_1}$ is $\Sigma^1_1$ (this fact, along with all the other basic facts about admissible ordinals, can be found in the book~\cite{chong2015}).

There are several properties of $F_{\omega_1}$ which make it interesting from the perspective of descriptive set theory. First, it is an analytic equivalence relation, all of whose equivalence classes are Borel and which has exactly $\aleph_1$-many equivalence classes (and which does not have perfectly many non-equivalent elements). Thus $F_{\omega_1}$ is a counterexample to a particularly strong generalization of Vaught's Conjecture. Second, Marker showed that $F_{\omega_1}$ is not induced by any continuous Polish group action~\cite{marker1988analytic} (later strengthened by Becker to show that $F_{\omega_1}$ is not induced by any Borel Polish group action~\cite{becker1994topological}). Thus $F_{\omega_1}$ is not a counterexample to the topological Vaught Conjecture. Further results on $F_{\omega_1}$ have been obtained by Montalb\'an~\cite{montalban2017degree}, Chan~\cite{chan2017countable} and others.

Our main result about $\Fw$ is the following.

\begin{theorem}
\label{thm:iff_nonconstructible}
The equivalence relation $\Fw$ is Borel graphable if and only if there is a non-constructible real.
\end{theorem}

Note that the forward direction of this theorem (as well as a couple of other results in this section) was first obtained by the first author in his PhD thesis~\cite{arant2019effective}.

Mansfield \cite{mansfield1975} established that
the existence of a non-constructible real is equivalent
to the non-existence of a $\Sigma^1_2$ wellordering of the reals,
which yields the following corollary to Theorem~\ref{thm:iff_nonconstructible}.

\begin{corollary}
The equivalence relation $\Fw$ is Borel graphable if and only if there is no $\Sigma^1_2$ wellordering of the reals.
\end{corollary}

As a corollary of the proof of this theorem, we will obtain a non-trivial example of an analytic equivalence relation which is \emph{provably} not Borel graphable. In particular, consider the equivalence relation $\Fwrel$ obtained by relativizing $F_{\omega_1}$ to each real and taking the direct sum of all the resulting equivalence relations. More precisely, let $\Fwrel$ denote the equivalence relation on $\Cantor \times \Cantor$ defined by
\[
(x, y) F^\text{rel}_{\omega_1}(x', y') \iff
(x=x' \ \text{ and } \ \omega_1^{x\oplus y}=\omega_1^{x\oplus y'}).
\]
We will show that, provably in $\ZFC$, $\Fwrel$ is not Borel graphable (see Corollary~\ref{cor:rel_not_graphable} below).

Theorem~\ref{thm:iff_nonconstructible} also has some 
bearing on the complexity of Borel graphability.

\begin{corollary}
There is no $\Pi^1_3$ formula which is provably equivalent to Borel graphability. I.e.\ there is no $\Pi^1_3$ formula $\varphi$ such that, provably in $\ZFC$, for every code $x$ for an analytic equivalence relation $E_x$, Borel graphability of $E_x$ is equivalent to $\varphi(x)$.
\end{corollary}

\begin{proof}
If there was such a formula $\varphi$, then by Shoenfield absoluteness, Borel graphability of any fixed analytic equivalence relation would be preserved under taking inner models. However, this is false: $F_{\omega_1}$ is Borel graphable if there is a non-constructible real but is not Borel graphable in the inner model $L$.
\end{proof}

It is tempting to think that we can strengthen the corollary above to simply state that Borel graphability is not expressible by a $\Pi^1_3$ formula. However, this does not obviously follow from Theorem~\ref{thm:iff_nonconstructible}. The problem is that it's possible that Borel graphability is provably $\Pi^1_3$, but exactly \emph{which} $\Pi^1_3$ formula it's equivalent to is independent of $\ZFC$. For example, it could be that Borel graphability is equivalent to one $\Pi^1_3$ formula when every real is constructible and equivalent to a different $\Pi^1_3$ formula when there is a non-constructible real. So the best we can conclude is the corollary as stated above: there is no single $\Pi^1_3$ formula which always works.

Although we cannot so far rule out the possibility that Borel graphability is $\Pi^1_3$, we suspect it is not. This motivates the following question.

\begin{question}
Is Borel graphability $\Sigma^1_3$-complete? That is, is the property that $x$ codes an analytic equivalence relation which is Borel graphable $\Sigma^1_3$-complete?
\end{question}

\subsection{When $\Fw$ is not Borel graphable}

We will begin by proving that when all reals are constructible, $\Fw$ is not Borel graphable (which, as mentioned above, was first proved in \cite{arant2019effective}). In order to prove this, we will first prove (in $\ZFC$) that $\Fw$ is not $\Delta^1_1$-graphable. By relativizing the proof of this fact and by combining it with a fact about constructible reals, we will obtain the full result. Along the way, we will prove the non-Borel graphability of $\Fwrel$.

In order to prove the non-$\Delta^1_1$ graphability of $\Fw$, we will need to review some definitions and theorems from computability theory and effective descriptive set theory. For reals $x, y \in \Cantor$, $x$ is \term{hyperarithmetically reducible} to $y$, denoted $x \leq_h y$, if $x$ is a $\Delta^1_1(y)$ subset of $\omega$. A real $x$ is simply called \term{hyperarithmetic} if it is hyperarithmetically reducible to $0$. The following standard theorem ensures that hyperarithmetic reducibility behaves like other computability-theoretic reducibility notions. For proofs of the statements in this theorem and a more thorough discussion of hyperarithmetic reducibility, see Section 2.4 of~\cite{chong2015} or Section II.5 of~\cite{sacks1990higher}.

\begin{theorem}
Hyperarithmetic reducibility has the following properties.
\begin{enumerate}
    \item Hyperarithmetic reducibility is transitive.
    \item If $x \leq_h z$ and $y \leq_h z$ then $x\oplus y \leq_h z$.
    \item For each $y$, there are countably many $x$ such that $x \leq_h y$. In particular, there are only countably many hyperarithmetic reals.
    \item If a set $A \subseteq \Cantor$ is $\Delta^1_1(x)$ and $x \leq_h y$ then $A$ is $\Delta^1_1(y)$. In particular, every set which is $\Delta^1_1$ relative to a hyperarithmetic parameter is $\Delta^1_1$.
\end{enumerate}
\end{theorem}

Given a real $x$, the \term{hyperjump} of $x$, denoted $\O^x$, is a complete $\Pi^1_1(x)$ subset of $\omega$, typically identified with the set of indices for well-orders computable from $x$. $\O^0$ is usually denoted simply as $\O$ and referred to as \term{Kleene's $\O$}. The hyperjump of a real $x$ is always strictly above $x$ for hyperarithmetic reducibility. The following standard theorem lists some useful properties of the function $x \mapsto \omega_1^x$, including its relationship to hyperarithmetic reducibility and the hyperjump. For a proof of the statements in this theorem, see Sections 2.4 and 2.5 of~\cite{chong2015} and II.5 and II.7 of~\cite{sacks1990higher}.

\begin{theorem}
\label{thm:omega1}
The function $x \mapsto \omega_1^x$ has the following properties.
\begin{enumerate}
    \item If $x \leq_h y$ then $\omega_1^x \leq \omega_1^y$. In particular, for every hyperarithmetic $x$ we have $\omega_1^x = \omega_1^{CK}$.
    \item If $\O^x \leq_h y$ then $\omega_1^x < \omega_1^y$. In particular, if $\O \leq_h x$ then $\omega_1^{CK} < \omega_1^x$.
    \item For every $x$, there are uncountably many $y$ such that $\omega_1^y = \omega_1^x$. In particular, there are uncountably many $y$ such that $\omega_1^y = \omega_1^{CK}$.
\end{enumerate}
\end{theorem}

We will also make use of the following two theorems. The second, sometimes known as Friedman's Conjecture, was originally proved by Martin~\cite{martin1976proof} and independently by Friedman (unpublished).

\begin{theorem}[Effective Perfect Set Theorem \cite{harrison1967}]
For every $\Sigma^1_1$ set $A \subseteq \Cantor$, either $A$ is uncountable or all elements of $A$ are hyperarithmetic.
\end{theorem}

\begin{theorem}[Friedman's Conjecture~\cite{martin1976proof}]
If $A \subseteq \Cantor$ is an uncountable $\Delta^1_1$ set then for every $x \geq_h \O$, there is some $y \in A$ such that $x \equiv_h y$.
\end{theorem}

Our use of the effective perfect set theorem and Friedman's Conjecture is encapsulated in the following corollary.

\begin{corollary}
\label{cor:friedman}
If $A \subseteq \Cantor$ is $\Delta^1_1$, then either every element of $A$ is hyperarithmetic or there is some element $x \in A$ such that $\omega_1^x > \omega_1^{CK}$.
\end{corollary}

\begin{proof}
By the effective perfect set theorem, either every element of $A$ is hyperarithmetic or $A$ is uncountable. In the latter case, we can apply Friedman's Conjecture to find some $x \in A$ such that $\O \leq_h x$. Theorem~\ref{thm:omega1} then implies that $\omega_1^{CK} < \omega_1^x$.
\end{proof}

We are now ready to prove that $\Fw$ is not $\Delta^1_1$-graphable.

\begin{proposition}
\label{prop:not_delta11}
$\Fw$ is not $\Delta^1_1$-graphable.
\end{proposition}

\begin{proof}
Suppose for contradiction that $G$ is a $\Delta^1_1$ graphing for $\Fw$.  Fix any hyperarithmetic real $a$. We claim that every element of the connected component of $a$ is hyperarithmetic. This is enough to finish the proof, since the class of hyperarithmetic reals is countable, but the $\Fw$-equivalence class of $a$ is uncountable by Theorem~\ref{thm:omega1}.

Suppose that there is a path in $G$ from $a$ to some real $b$. We must show that $b$ is hyperarithmetic. We prove this by induction on the length of the path. The base case---when the length of the path is $0$, i.e.\ when $a = b$---follows from our choice of $a$. For the inductive case, let $c$ be the element of the path immediately preceding $b$. We know by the inductive assumption that $c$ is hyperarithmetic. Thus the set of neighbors of $c$ is $\Delta^1_1(c) = \Delta^1_1$. So by Corollary~\ref{cor:friedman} to Friedman's Conjecture, either every neighbor of $c$ is hyperarithmetic (so, in particular, $b$ is hyperarithmetic) or there is some neighbor $x$ of $c$ such that $\omega_1^x > \omega_1^{CK}$. But this is impossible: since $G$ is a graphing of $\Fw$, $x$ must be $\Fw$-equivalent to $c$ and thus must satisfy $\omega_1^x = \omega_1^c = \omega_1^{CK}$.
\end{proof}

\begin{corollary}
\label{cor:rel_not_graphable}
$\Fwrel$ is not Borel graphable.
\end{corollary}

\begin{proof}
If $F^{\text{rel}}_{\omega_1}$ is Borel graphable by some $\Delta^1_1(x)$ graph $G$, then one can obtain a contradiction from the relativized version of the above argument applied to the $G$-connected component of $(x, x)$.  
\end{proof}

It is perhaps worth examining why a relativization of the proof of Proposition~\ref{prop:not_delta11} does not show that $\Fw$ is not Borel graphable. Suppose that $\Fw$ was $\Delta^1_1(x)$-graphable for some real $x$. The natural way to try to relativize the proof of Proposition~\ref{prop:not_delta11} is to try to prove that every element of the connected component of $x$ is hyperarithmetically reducible to $x$. To see what goes wrong, let's consider the set of neighbors of $x$. This set is $\Delta^1_1(x)$ and thus by the (relativized) effective perfect set theorem, either all its elements are hyperarithmetically reducible to $x$ or it is uncountable. We would like to apply Friedman's Conjecture to rule out the latter case, but there is a problem. The relativized version of Friedman's Conjecture applied in this context gives us a neighbor $y$ of $x$ such that $y \oplus x \geq_h \O^x$ and thus $\omega_1^{y\oplus x} > \omega_1^x$. However, unlike in the unrelativized proof, this does not immediately yield a contradiction. The problem is that $\omega_1^{y \oplus x}$ can in general be much larger than $\omega_1^y$ and to get a contradiction we need that $\omega_1^y > \omega_1^x$ (which would contradict the $\Fw$-equivalence of $x$ and $y$).

In general, this problem cannot be overcome: we have already mentioned that it is consistent with $\ZFC$ that $\Fw$ is Borel graphable. However, when all reals are constructible, we can get around this problem by invoking the following theorem about hyperarithmetic reducibility of reals in $L$.

\begin{theorem}[Guaspari, Kechris, Sacks]
\label{thm:constructible_reals}
For every real $x \in L$, there is some real $y \in L$ such that $x \leq_h y$ and $y$ is of the least degree of hyperarithmetic reducibility in its $\Fw$-equivalence class. In other words, for all reals $z$, if $\omega_1^y = \omega_1^z$ then $y \leq_h z$.\footnote{
This theorem can be deduced, for example, from (a relativized version of)
Theorem 3.6.8 and Proposition 4.3.4 in \cite{chong2015}.
}
\end{theorem}

\begin{proposition}
Suppose that all reals are constructible. Then $\Fw$ is not Borel graphable.
\end{proposition}

\begin{proof}
The proof closely follows the proof of Proposition~\ref{prop:not_delta11}. Suppose for contradiction that $\Fw$ is $\Delta^1_1(x)$-graphable for some real $x$. By Theorem~\ref{thm:constructible_reals}, we may assume that if $\omega_1^x = \omega_1^y$ then $x \leq_h y$.

We claim that every element of the connected component of $x$ is hyperarithmetically reducible to $x$. As before, this is enough to finish the proof: there are only countably many reals hyperarithmetically reducible to $x$, but the $\Fw$-equivalence class of $x$ is uncountable.

By induction, it is enough to check that if $y$ is both $\Fw$-equivalent and hyperarithmetically reducible to $x$ then all the neighbors of $y$ are also hyperarithmetically reducible to $x$ (it follows immediately that all the neighbors are $\Fw$-equivalent to $x$ because they are all in the same connected component in a graphing of $\Fw$). Note that the set of neighbors of $y$ is a $\Delta^1_1(x\oplus y) = \Delta^1_1(x)$ set. Thus by the relativized effective perfect set theorem, either every neighbor of $y$ is hyperarithmetically reducible to $x$ or $y$ has uncountably many neighbors. 

We now claim the latter case is impossible, which is enough to finish the proof. For contradiction, assume the latter case holds---i.e.\ $y$ has uncountably many neighbors. We can apply the relativized version of Friedman's Conjecture to obtain a neighbor $z$ of $y$ such that $z \oplus x \geq_h \O^x$ and thus $\omega_1^{z\oplus x} > \omega_1^x$. Note that $z$ and $x$ are $\Fw$-equivalent and thus $\omega_1^x = \omega_1^z$. By our choice of $x$, this implies that $x \leq_h z$ and thus $z\oplus x \equiv_h z$. Therefore $\omega_1^z = \omega_1^{z\oplus x} > \omega_1^x$, a contradiction of the $\Fw$-equivalence of $x$ and $z$.
\end{proof}

\subsection{When $\Fw$ is Borel graphable}

We will now prove that if there is a non-constructible real then $\Fw$ has the Borel witness coding property and thus is Borel graphable. Instead of directly showing that $\Fw$ has the Borel witness coding property if there is a non-constructible real, we will first prove that $\Fw$ has the Borel witness coding property under a somewhat stronger assumption and then use this to prove the full result via an absoluteness argument.\footnote{We are grateful to Gabe Goldberg for pointing out this absoluteness argument.}

In our proof, we will make use of a technique known as Kumabe-Slaman forcing. We will now briefly review the main facts about Kumabe-Slaman forcing that we will use in our proof; a more complete introduction to Kumabe-Slaman forcing, including proofs of its key properties, can be found in the appendix.

Kumabe-Slaman forcing is a type of forcing which produces a real. Its key property is as follows: given a countable transitive model $M$ of $\ZFC$ (or of a sufficiently large finite fragment of $\ZFC$) and a real $x$ not in $M$, there is some $g$ which is Kumabe-Slaman generic over $M$ such that $g\oplus x$ encodes essentially any information, including information about $g$ itself. This is made more precise in Propositions~\ref{prop:encode1} and~\ref{prop:encode2} in the appendix, but essentially boils down to the following. Suppose that $\phi_0, \phi_1, \ldots$ is a sequence of formulas in the forcing language for Kumabe-Slaman forcing over $M$. Then there is some $g$ which is Kumabe-Slaman generic over $M$ such that $g\oplus x$ computes the truth of each $\phi_n$ in $M[g]$, uniformly in $n$.

A second key fact about Kumabe-Slaman forcing is that it preserves admissibility. More precisely, suppose $M$ is a transitive model of $\ZFC$ and $x$ is a real in $M$. Then for every $g$ which is Kumabe-Slaman generic over $M$, $\omega_1^{x\oplus g} = \omega_1^x$. See Theorem~\ref{thm:admissible} in the appendix for a proof of this fact.

\begin{lemma}
\label{lemma:omega_1_L}
Suppose that there is a non-constructible real and $\omega_1^L = \omega_1$. Then $\Fw$ has the Borel witness coding property.
\end{lemma}

\begin{proof}
Let $a$ be a non-constructible real. We will show that for all reals $x$ and $b$, there is some real $y$ such that $\omega_1^x = \omega_1^y$ and $y\oplus a$ computes both $b$ and a witness to the $\Fw$-equivalence of $x$ and $y$. Thus $\Fw$ has the Borel witness coding property via the function $x \mapsto x\oplus a$.

Fix reals $x$ and $b$ as discussed above. Since $\omega_1^L = \omega_1$, $\omega_1^x$ is countable in $L$. Thus we can find some real $y_0 \in L$ such that $\omega_1^x = \omega_1^{y_0}$ (this follows from the fact that for a countable ordinal $\alpha$, the property ``there is some $x$ such that $\omega_1^x = \alpha$'' is $\Sigma^1_1$ relative to a code for $\alpha$ and thus if $\alpha$ is countable in $L$, this property is absolute between $L$ and $V$).

Let $\alpha$ be some countable ordinal large enough that $y_0 \in L_\alpha$ and $L_\alpha \models \ZFC$ (or at least a sufficiently large finite fragment of $\ZFC$). Form $y_1$ using Kumabe-Slaman forcing over $L_\alpha$ in such a way that $y_1\oplus a$ computes both $b$ and some additional information that we will explain below (to see why we can at least choose $y_1$ so that $y_1\oplus a$ computes $b$, refer to Proposition~\ref{prop:encode1}).

Now define $y = y_0\oplus y_1$. We claim that $y$ has the desired properties. First, since Kumabe-Slaman forcing preserves admissibility, we have
\[
    \omega_1^{y_0\oplus y_1} = \omega_1^{y_0} = \omega_1^x
\]
and thus $x$ and $y$ are $\Fw$-equivalent. Second, by our choice of $y_1$, $y\oplus a$ computes $b$. It only remains to check that $y \oplus a$ computes a witness to the $\Fw$-equivalence of $x$ and $y$. Actually, since we were a bit vague above about how $y_1$ is chosen, what we really need to do is explain why we can find some $y_1$ which is Kumabe-Slaman generic over $L_\alpha$ and for which $y_1\oplus a$ can compute such a witness.

The most straightforward way to check that we can find $y_1$ for which $y_1\oplus a$ computes the desired witness is to write down exactly what such a witness consists of and check that each piece of information needed can be expressed in the forcing language for Kumabe-Slaman forcing over $L_\alpha$ (or can be replaced with a piece of information which can be so expressed) and is absolute between $L_\alpha[y_1]$ and $V$. However, this is fairly tedious to do, so we will use a slightly cheap trick. The key point is that it is enough for $y_1$ to encode both a witness to the $\Fw$-equivalence of $x$ and $y_0$ and a witness to the $\Fw$-equivalence of $y_0$ and $y$. The former does not depend at all on $y_1$ and so it can be encoded using the method of Proposition~\ref{prop:encode1}. Thus we will just explain how to encode the latter.

For this, fix a $\Pi^0_1$ set $C \subseteq \Cantor\times\Cantor\times\Baire$ such that $C$ is an unfolding of $F_{\omega_1}$ (which is possible since $\Fw$ is $\Sigma^1_1$, as mentioned above). For any $u, v \in \Cantor$, let $\proj(C, u, v)$ denote the set $\{z \mid (u, v, z) \in C\}$ and recall that, formally, a witness to the $\Fw$-equivalence of $u$ and $v$ is exactly an element of $\proj(C, u, v)$. Note that whenever $u$ and $v$ really are $\Fw$-equivalent, the set $\proj(C, u, v)$ is a closed, nonempty subset of $\Baire$ and thus has a leftmost element. Let $\phi$ be a formula expressing the values of the digits of this leftmost element---i.e.\ for each $n, m \in \N$ and $u,v \in \Cantor$, $\phi(n, m, u, v)$ states that the $m$ is the $n^\text{th}$ digit of the leftmost element of $\proj(C, u, v)$ (if it exists). Note that since $C$ is $\Pi^0_1$, $\phi$ can be written as the conjunction of a $\Sigma^1_1$ formula and a $\Pi^1_1$ formula.

Now observe two things. First $y_0$ and $y$ really are $\Fw$-equivalent and so $\proj(C, y_0, y)$ has a leftmost element. Second, because $L_\alpha[y_1]$ is a well-founded model of $\ZFC$, it is correct about the digits of this leftmost element (since $\phi$ is $\Sigma^1_1\land \Pi^1_1$). Thus using Proposition~\ref{prop:encode2}, we can choose $y_1$ so that all the digits of this leftmost element are encoded into $y_1\oplus a$.
\end{proof}

\begin{proposition}
\label{prop:graphable}
If there is a non-constructible real then $\Fw$ has the Borel witness coding property and hence is Borel graphable with diameter $2$.
\end{proposition}

\begin{proof}
By Lemma~\ref{lemma:omega_1_L} above, we may assume that $\omega_1^L < \omega_1$. Thus $(\Cantor)^L$ is countable and so there is a real $g$ which is Cohen generic over $L$. Since Cohen generics preserve cardinals, $L[g] \models \omega_1^L = \omega_1$ and thus $L[g]$ satisfies the conditions required by Lemma~\ref{lemma:omega_1_L}. In other words, $\Fw$ has the Borel witness coding property in $L[g]$. Since having the Borel witness coding property can be expressed by a $\Sigma^1_3$ formula, Shoenfield absoluteness implies that it is preserved upwards from $L[g]$ to $V$. Thus $\Fw$ has the Borel witness coding property.
\end{proof}


\section{Isomorphism of countable structures}
\label{sec:isomorphism}

Theorem~\ref{thm:loiso} from Section~\ref{sec:basics} states that the isomorphism relation on the space of countable linear orders is Borel graphable. This result begs for generalization---in particular, is there anything special about linear orders here or does it hold for the isomorphism relation on other classes of structures? If so, which ones? In this section, we will answer these questions by proving the following theorem.

\begin{notation}
Given a countable language $\lang$, $\Mod(\lang)$ denotes the Polish space of $\lang$-structures with domain $\N$. See Section 3.6 of the book~\cite{gao2009invariant} of Gao for a more detailed definition of this space and a proof that it is Polish.
\end{notation}

\begin{theorem}
\label{thm:isomorphism}
For any countable language $\lang$, the isomorphism relation on $\Mod(\lang)$ is Borel graphable with diameter $2$.
\end{theorem}

As a consequence of this theorem, we also obtain Borel graphability of piecewise isomorphism and bi-embeddability for structures in an arbitrary countable language, thus generalizing Corollary~\ref{cor:lo}. Also, in Section~\ref{sec:graphic} we will restate this theorem in the language of Polish group actions and investigate to what extent it can be generalized in that setting.

\begin{remark}
After we wrote this paper, it was pointed out to us that the proof of Theorem~\ref{thm:isomorphism} has substantial overlap with the proof of a theorem of Julia Knight. In particular, Knight~\cite{knight1986degrees} (generalizing a theorem of Marker~\cite{marker1982degrees}) showed that if $\lang$ is a countable language and $M \in \Mod(\lang)$ then either $M$ is trivial in a certain sense or for every $x \in \Cantor$ which computes $M$, there is a structure $N \in \Mod(\lang)$ which is isomorphic to $M$ and Turing equivalent to $x$. The bulk of our proof is devoted to proving something very similar, but slightly stronger: we will show that in the case where $M$ is non-trivial, we can actually find $N$ which is isomorphic to $M$ and which computes both $x$ and an isomorphism $M \to N$. Furthermore, our proof of this fact proceeds along very similar lines to those of Knight's proof, with the addition of essentially one idea: the ``queue of information'' described in the proof of Lemma~\ref{lemma:non_trivial_isomorphism} below.
\end{remark}

\begin{remark}
It might seem at first glance that Theorem~\ref{thm:isomorphism} is not quite a generalization of the Borel graphability of $\loiso$. After all, the underlying space $\lo$ does not consist of all countable structures in the language of linear orders, but just those that satisfy the theory of linear orders. However, note that for any first order $\lang$-theory $T$, the set of models of $T$ is an isomorphism invariant subset of $\Mod(\lang)$. Thus the isomorphism relation restricted to this set is Borel graphable as well and so Theorem~\ref{thm:isomorphism} directly implies Borel graphability of $\loiso$.
\end{remark}

\begin{remark}
Note that any countable language $\lang$ can be represented by an element of $\Cantor$ and thus we may talk about computability of $\lang$ and of models of $\lang$ with domain $\N$.
\end{remark}

\subsection{Isomorphism is Borel graphable}

To prove Theorem~\ref{thm:isomorphism} we will partition $\Mod(\lang)$ into two isomorphism invariant Borel subsets such that the isomorphism relation is essentially trivial when restricted to the first subset and has the Borel witness coding property when restricted to the second subset. Roughly speaking, the first part of this partition consists of those $\lang$-structures for which, aside from a finite set, all elements look identical. In order to make this precise, it will be helpful to introduce some terminology.

Suppose $M$ is an $\lang$-structure, $F$ is a finite subset of $M$ and $\bar{a}$ is a tuple of distinct elements of $M$. The \term{quantifier free type of $\bar{a}$ over $F$} is the set of quantifier-free formulas $\phi(\bar{x})$ with parameters from $F$ such that $M \models \phi(\bar{a})$. We will call an $\lang$-structure $M$ \term{trivial} if there is some finite set $F \subset M$ such that for all tuples $\bar{a}$ and $\bar{b}$ of distinct elements of $M \setminus F$, if $\bar{a}$ and $\bar{b}$ have the same length then they have the same quantifier-free type over $F$.

\begin{remark}
It is possible to show that an $\lang$-structure is trivial if and only if its isomorphism class is countable, but we will not directly use this characterization of triviality.
\end{remark}

It is easy to see that the set of trivial $\lang$-structures is both isomorphism invariant and Borel. Thus to prove that the isomorphism relation on $\Mod(\lang)$ is Borel graphable, it is enough to prove separately that it is graphable when restricted to the trivial $\lang$-structures and when restricted to the non-trivial $\lang$-structures (because we can then just take the disjoint union of the two resulting graphs). The following two lemmas accomplish this.

\begin{lemma}
\label{lemma:trivial_isomorphism1}
If $M$ and $N$ are trivial $\lang$-structures which are isomorphic then there is a computable isomorphism between them. Hence isomorphism of trivial $\lang$-structures is Borel.
\end{lemma}

\begin{proof}
Let $f \colon M \to N$ be an isomorphism and let $F \subset M$ be a finite set witnessing that $M$ is trivial. We claim that any bijection $g \colon M \to N$ which agrees with $f$ on $F$ is an isomorphism. Given this claim, we can easily finish the proof: the claim implies that to compute an isomorphism from $M$ to $N$, we just need to know where to send the elements of $F$ and since $F$ is finite, this is just a finite amount of information.

We now prove the claim. Fix tuples $\bar{a}$ and $\bar{b}$ of distinct elements of $F$ and $M\setminus F$, respectively, and let $\phi(\bar{x}, \bar{y})$ be a quantifier-free $\lang$-formula. It is enough to show that $M \models \phi(\bar{a}, \bar{b})$ if and only if $N \models \phi(g(\bar{a}), g(\bar{b}))$. Let $\bar{c}$ be a tuple in $M$ such that $f(\bar{c}) = g(\bar{b})$ and note that since $f$ and $g$ are injective and agree on $F$, $\bar{c}$ is disjoint from $F$. We can now calculate
\begin{align*}
  M \models \phi(\bar{a}, \bar{b}) &\iff M \models \phi(\bar{a}, \bar{c}) &\text{(triviality of $M$)}\\
                                   &\iff N \models \phi(f(\bar{a}), f(\bar{c})) &\text{($f$ is an isomorphism)}\\
                                   &\iff N \models \phi(g(\bar{a}), g(\bar{b})) &\text{($f$ and $g$ agree on $F$ and $f(\bar{c}) = g(\bar{b})$)}
\end{align*}
which finishes the proof.
\end{proof}

\begin{lemma}
\label{lemma:non_trivial_isomorphism}
The isomorphism relation restricted to non-trivial $\lang$-structures has the Borel witness coding property.
\end{lemma}

\begin{proof}
Fix a non-trivial $\lang$-structure $M$ and a real $r \in \Cantor$. We will build an $\lang$-structure $N$ which is isomorphic to $M$ such that the Turing jump of $N \oplus \lang$ computes both $r$ and an isomorphism between $M$ and $N$. This is enough to finish the proof, since it shows that the isomorphism relation on non-trivial $\lang$-structures has the Borel witness coding property via the function which sends a structure $M$ to $(M\oplus \lang)'$.

When describing our construction of $N$, it will be helpful to think of the elements of $M$ not as natural numbers (which they technically are) but as elements of some abstract set. Also, to follow the construction, it helps to keep in mind the following mental image: we will build $N$ by picking up the elements of $M$ one at a time and placing them on natural numbers. We will begin by placing some element of $M$ on $0$, then placing some element on $1$ and so on, until every natural number has some element of $M$ placed on it. We can then define $N$ as the $\lang$-structure on $\N$ where each natural number behaves like the element of $M$ which has been placed on it. To ensure that $N$ is isomorphic to $M$, we just need to ensure that every element of $M$ is eventually placed somewhere.

We will carry out this construction in a series of stages, on each of which we will place elements of $M$ on a finite number of natural numbers. On each stage, we will first place a few elements of $M$ to code one bit of information and then place one more element of $M$ to help ensure that every element is eventually placed somewhere.

Here's what this means more formally. We will build a bijection $f \colon \N \to M$ and define $N$ by setting, for each relation $R$ in $\lang$,
\[
  R^N(n_0, \ldots, n_k) \iff R^M(f(n_0), \ldots, f(n_k)).
\]
Placing an element $a$ of $M$ on a number $n$ corresponds to setting $f(n) = a$ and ensuring that every element of $M$ is eventually placed somewhere corresponds to ensuring that $f$ is surjective. We will build $f$ in a series of stages, where at the end of each stage we will have determined the values of $f$ on some finite initial segment of $\N$.

\medskip\noindent\textit{The construction.}
We will now describe in more detail what happens on each stage of the construction. Suppose that we are on stage $n$ and that so far, we have placed elements of $M$ on the numbers $0, 1, 2, \ldots, k$ (i.e.\ that we have defined $f$ on these numbers). Let $\bar{a} = (a_0, \ldots, a_k)$ denote the elements of $M$ which were placed on these numbers and set $F = \{a_0, \ldots, a_k\}$. The key idea of the construction is that we can code information using our choice of the quantifier-free type over $F$ of the next few elements that we place.

Since $M$ is not trivial, there are two tuples of elements of $M\setminus F$ which have different quantifier-free types over $F$. Let $\phi(\bar{x}, \bar{y})$ be the first quantifier free formula (according to some reasonable ordering) such that there are tuples $\bar{b}$ and $\bar{c}$ in $M\setminus F$ for which $M \models \phi(\bar{a}, \bar{b})$ and $M \notmodels \phi(\bar{a}, \bar{c})$. If the next bit of information to be coded is $0$ then place the elements of $\bar{b}$ on the next available numbers (i.e.\ on the numbers $k + 1, \ldots, k + l$, where $l$ is the length of the tuple $\bar{b}$). Otherwise, place the elements of $\bar{c}$ on these numbers. We will explain below exactly what information should be encoded on each stage.

To finish the stage, let $d$ be the least element of $M$ not yet placed anywhere and place it on the next available number (i.e.\ on the number $k + l + 1$).

\medskip\noindent\textit{Verification.}
We will now show that one jump of $N \oplus \lang$ can compute both $r$ and an isomorphism between $M$ and $N$. To do so, we will first show that one jump of $N \oplus \lang$ can decode all the information that was encoded during the construction and then explain what sort of information needs to be encoded on each stage. The main idea of the proof is just that it is not very hard to use $N$ to unravel what happened on each stage of the construction and thus to determine what information was encoded.

Suppose that, using $(N \oplus \lang)'$, we have so far determined what happened on the first $n$ stages of the construction. This means, in particular, that we have determined both what information was encoded on those first $n$ stages and the last number, $k$, to have had some element of $M$ placed on it by the end of stage $n$. To determine what happened on stage $n + 1$, we can use $(N \oplus \lang)'$ to search for the first quantifier-free formula $\phi$ (with parameters consisting of the first $k$ elements of $N$) for which two tuples of elements in $N$ give different values. The number $l$ of free variables in this tuple tells us how many elements were placed on stage $n + 1$ (namely, $l + 1$) and the truth value of $\phi$ on the next $l$ elements of $N$ tells us whether a $0$ or a $1$ was encoded on stage $n + 1$.

We will now explain what information to encode on each stage. Of course, we need to encode the bits of the real $r$. We also need to encode an isomorphism between $M$ and $N$, but the values of such an isomorphism may depend on what elements of $M$ we choose to place where. We can solve this apparent circular dependence by maintaining a ``queue'' of information. On each stage, we will pull one bit of information off the front of this queue and encode it and then add some more bits of information to the back of the queue to encode during future stages.\footnote{More precisely, we will maintain a list of bits $b_0, b_1, \ldots, b_k$ and a counter $m$. On each stage, we will encode bit $b_m$, increment $m$ by one (i.e.\ replace $m$ with $m + 1$) and then extend the list $b_0, \ldots, b_m$ by adding a finite number of extra bits to the end of it. In this way, we will ensure that all bits are eventually encoded.}

During each stage, we determine the placement of finitely many elements of $M$. At the end of the stage, we can add information encoding these placements into the queue, along with the next bit of $r$  (perhaps using some coding scheme so that we can tell when a bit of information is a bit of $r$ rather than information about the placement of elements).
\end{proof}

\subsection{Piecewise isomorphism and bi-embeddability are Borel graphable}

We will now show that for an arbitrary countable relational language $\lang$, the equivalence relations of piecewise isomorphism and bi-embeddability on $\Mod(\lang)$ are both Borel graphable.\footnote{For the definition of piecewise isomorphism, see the introduction. Note that the requirement that $\lang$ is relational is necessary for piecewise isomorphism since the definition of piecewise isomorphism only really makes sense for relational languages. It is not needed for bi-embeddability.} However, the Borel graphability of these two equivalence relations does not follow quite as easily from the proof of graphability of isomorphism as it did in the case of linear orders in Section~\ref{sec:basics}. Here's why.

To prove Borel graphability of isomorphism, we divided $\Mod(\lang)$ into the set of trivial $\lang$-structures and the set of non-trivial $\lang$-structures. We then showed that isomorphism of trivial $\lang$-structures is Borel and isomorphism of non-trivial $\lang$-structures has the Borel witness coding property. The first problem is that while the latter result implies that piecewise isomorphism (or bi-embeddability) restricted to the non-trivial $\lang$-structures has the Borel witness coding property, the former result does not immediately imply that piecewise isomorphism (or bi-embeddability) is Borel when restricted to the trivial $\lang$-structures (since it is possible that a superequivalence relation of a Borel equivalence relation is not itself Borel). The second problem is that it is typically not the case that the set of trivial $\lang$-structures is invariant for piecewise isomorphism or bi-embeddability.

To fix the first problem, we will prove directly that piecewise isomorphism and bi-embeddability \emph{are} Borel when restricted to the trivial $\lang$-structures. We will then see that the second problem does not present a major obstacle.

\begin{lemma}
\label{lemma:trivial_isomorphism2}
Suppose that $\lang$ is a relational language and that $M$ and $N$ are trivial $\lang$-structures.
\begin{enumerate}
\item If $M$ and $N$ are piecewise isomorphic, then there is a computable piecewise isomorphism.
\item If $M$ and $N$ are bi-embeddable, then there is a computable bi-embedding.
\end{enumerate}
\end{lemma}

\begin{proof}
In both cases, we can imitate the proof of Lemma~\ref{lemma:trivial_isomorphism1} with a few extra details added.

We will begin with (1). Suppose $M$ and $N$ are piecewise isomorphic, as witnessed by $(f, A_1, \ldots, A_n)$, and trivial, as witnessed by $F \subset M$ and $E \subset N$. By replacing $F$ with $F \cup f^{-1}(E)$, we may assume that $E \subseteq f(F)$.

Now consider the partition $B_1, \ldots, B_m$ of $M$ which puts each element of $F$ into a part by itself and has one more part consisting of $M \setminus F$. Let $g \colon M \to N$ be any bijection which agrees with $f$ on $F$. We claim that $(g, B_1, \ldots, B_m)$ is a piecewise isomorphism of $M$ and $N$. As in the proof of Lemma~\ref{lemma:trivial_isomorphism1}, this is enough to finish the proof: to compute such a $g$ (and the partition $B_1,\ldots,B_m$), we just need to remember the elements of $F$ and where each one is mapped by $f$, which is a finite amount of information.

Since $g$ agrees with $f$ on $F$ and each element of $F$ has been put into a piece of the partition by itself, it suffices to fix a quantifier free $\lang$-formula $\phi(\bar{x})$ and a tuple $\bar{a}$ of distinct elements of $M \setminus F$ and show that $M \models \phi(\bar{a})$ if and only if $N \models \phi(g(\bar{a}))$. Let $\bar{b}$ be any tuple of distinct elements of $M\setminus F$
all of which are in the same part of the partition $(A_1,\ldots,A_n)$ and
which has the same length as $\bar{a}$. Note that by our assumption on $F$, $f(\bar{b})$ and $g(\bar{a})$ are both contained in $N\setminus E$. We can now calculate
\begin{align*}
M \models \phi(\bar{a}) &\iff M \models \phi(\bar{b}) &\text{(triviality of $M$)}\\
                        &\iff N \models \phi(f(\bar{b})) &\text{($\bar{b}$ comes from a single piece of the partition)}\\
                        &\iff N \models \phi(g(\bar{a})) &\text{(triviality of $N$)}
\end{align*}
which finishes the proof of (1).

We will now handle (2). It is enough to show that for every embedding of one trivial $\lang$-structure into another, there is also a computable embedding. To that end, suppose $N$ is trivial, as witnessed by $E \subset N$, and $f \colon M \to N$ is an embedding. Define $F = f^{-1}(E)$. We claim that any injection $g \colon M \to N$ which agrees with $f$ on $F$ and which avoids $E \setminus \range(f)$ is an embedding.

Fix a quantifier free $\lang$-formula $\phi$ and tuples $\bar{a}$ in $F$ and $\bar{b}$ in $M\setminus F$. It suffices to show that $M \models \phi(\bar{a}, \bar{b})$ if and only if $N \models \phi(g(\bar{a}), g(\bar{b}))$. Note that by definition of $F$ and our assumption on $g$, both $f(\bar{b})$ and $g(\bar{b})$ are disjoint from $E$. We can now calculate
\begin{align*}
M \models \phi(\bar{a}, \bar{b}) &\iff N \models \phi(f(\bar{a}), f(\bar{b})) &\text{($f$ is an embedding)}\\
                        &\iff N \models \phi(f(\bar{a}), g(\bar{b})) &\text{(triviality of $N$)}\\
                        &\iff N \models \phi(g(\bar{a}), g(\bar{b})) &\text{($f$ and $g$ agree on $F$)}
\end{align*}
which finishes the proof of (2).
\end{proof}

\begin{corollary}
For any countable relational language $\lang$, the piecewise isomorphism and bi-embeddability equivalence relations on $\Mod(\lang)$ are both Borel graphable with diameter $2$.
\end{corollary}

\begin{proof}
We only give the proof for piecewise isomorphism since the proof for bi-embeddability is identical. Let $f\colon \Mod(\lang) \to \Cantor$ be a function witnessing that the isomorphism relation restricted to non-trivial $\lang$-structures has the Borel witness coding property. Define a graph $G$ on $\Mod(\lang)$ by putting an edge between $M$ and $N$ if either of the following holds:
\begin{enumerate}
\item There is a computable piecewise isomorphism between $M$ and $N$.
\item $f(M)$ computes an $\lang$-structure $M'$, an isomorphism between $M$ and $M'$ and a piecewise isomorphism between $M'$ and $N$ (or similarly, but with the roles of $M$ and $N$ reversed).
\end{enumerate}
Now suppose that $M$ and $N$ are piecewise isomorphic. We will show that there is a path of length at most $2$ in $G$ connecting $M$ to $N$.

\medskip\noindent\textit{Case 1: both $M$ and $N$ are trivial.} By Lemma~\ref{lemma:trivial_isomorphism2}, $M$ and $N$ are adjacent in $G$.

\medskip\noindent\textit{Case 2: at least one of $M$ and $N$ is not trivial.} Without loss of generality, assume $M$ is not trivial. Fix a piecewise isomorphism, $g$, of $M$ and $N$. By our choice of $f$ and the fact that $M$ is not trivial, there is some $M'$ isomorphic to $M$ such that $f(M')$ computes $g$, $M$ and an isomorphism between $M'$ and $M$. We claim that $M'$ is adjacent to both $M$ and $N$. To see why $M'$ is adjacent to $N$, note that $f(M')$ computes $M$, an isomorphism between $M'$ and $M$ and a piecewise isomorphism between $M$ and $N$. To see why $M'$ is adjacent to $M$, note that $f(M')$ computes an isomorphism between $M'$ and $M$ and a piecewise isomorphism between $M$ and itself (e.g.\ the identity).
\end{proof}


\section{Graphic Polish groups}
\label{sec:graphic}

In this section, we will study the graphability of orbit equivalence relations arising from Polish group actions. In particular, we will consider the following class of Polish groups.

\begin{definition}
A Polish group $\Gamma$ is \term{graphic} if every Borel action of $\Gamma$ on a Polish space has a Borel graphable orbit equivalence relation.
\end{definition}

One motivation for this comes from the results of the previous section, which can be restated in the language of Polish group actions. In particular, let $S_\infty$ denote the Polish group of permutations of $\N$. For any language $\lang$, there is a natural action of $S_\infty$ on $\Mod(\lang)$, known as the \term{logic action of $S_\infty$ on $\Mod(\lang)$}, and the orbit equivalence relation of this action is exactly the isomorphism relation on $\Mod(\lang)$.

Becker and Kechris showed in~\cite{becker1993borel} that if $\lang$ is a language with relations of unbounded arity, then the logic action of $S_\infty$ on $\Mod(\lang)$ is universal for $S_\infty$. More precisely, for every Borel action of $S_\infty$ on a Polish space $X$, there is a Borel $S_\infty$-invariant embedding of $X$ into $\Mod(\lang)$. It follows that there is an invariant Borel reduction from the orbit equivalence relation on $X$ to the isomorphism relation on $\Mod(\lang)$. Thus Theorem~\ref{thm:isomorphism} (together with the fact that Borel graphability is preserved by invariant Borel reductions) implies that every Polish group action of $S_\infty$ has a Borel graphable orbit equivalence relation---in other words, that $S_\infty$ is graphic.

It is also not hard to see that $S_\infty$ is not the only graphic Polish group. For example, a Polish group is called \term{tame} if all of its Polish group actions induce Borel equivalence relations. It is known that there are many tame Polish groups, including all Polish groups which are countable or locally compact. Since any tame group is obviously graphic, this gives us many examples of graphic Polish groups ``for free.''

These observations raise an obvious question: are all Polish groups graphic? If so, this would constitute a far-reaching generalization of Theorem~\ref{thm:isomorphism}. So far, we have been unable to answer this question, but the results of this section can be read as a partial answer; in particular, we will show that many Polish groups are graphic and give a few conditions which are sufficient to imply that all Polish groups are graphic.

We now make two observations that are helpful when verifying that a Polish group $\Gamma$ is graphic. First, instead of checking all Borel actions of $\Gamma$, it is enough to just check the continuous actions of $\Gamma$. This follows from a result of Becker and Kechris~\cite{becker1993borel}, which states that if $\Gamma \acts X$ is a Borel action of $\Gamma$ on a Polish space $X$, then there is a Polish topology on $X$ which makes the action continuous and does not change the class of Borel subsets of $X$. Second, it is enough to just check actions of $\Gamma$ for which all orbits are uncountable. This follows from the following well-known (folklore) proposition.

\begin{proposition}
\label{prop:countable_orbit1}
If $\Gamma \acts X$ is a continuous action of a Polish group $\Gamma$ on a Polish space $X$, then
\[
X_C =\{x\colon \Gamma\cdot x  \ \textrm{is countable}\}
\]
is Borel and $E^X_\Gamma | X_C$ is Borel.

\end{proposition}

\begin{proof}
Fix a countable dense subset $\{g_n\}$ of $\Gamma$ and a countable basis $\{U_n\}$ of the identity of $\Gamma$. We claim that $x \in X_C$ if and only if the stabilizer $\Gamma_x$ of $x$ contains an open neighborhood of the identity. Granting this claim, we have
\[
x\in X_C \iff \exists m \forall k (g_k \in U_m \implies g_k\cdot x = x),
\]
and thus $X_C$ is Borel. Moreover, if $\Gamma_x$ contains a neighborhood of the identity then every coset of $\Gamma_x$ has nonempty interior. It follows that for $x,y \in X_C$,
\[
x E^X_\Gamma y \iff \exists k (g_k\cdot x = y),
\]
and so $E^X_\Gamma | X_C$ is Borel.

To see why $x \in X_C$ iff $\Gamma_x$ contains an open neighborhood of the identity, we can reason as follows. First note that $x \in X_C$ iff $\Gamma_x$ has countable index in $\Gamma$. If $\Gamma_x$ has countable index then it must be non-meager (otherwise all its cosets would be meager and thus $\Gamma$ would be covered by countably many meager sets). Thus by a theorem of Pettis (Theorem 9.9 of~\cite{kechris1995classical}), $\Gamma_x$ contains an open neighborhood of the identity. In the other direction, if $\Gamma_x$ contains an open neighborhood of the identity then all of its cosets have nonempty interior. Since $\Gamma$ is separable, it follows from this that there must only be countably many cosets of $\Gamma_x$.
\end{proof}

\subsection{Closure properties of graphic Polish groups}

We will begin by studying the closure properties of the class of graphic Polish groups. First, we note a few closure properties which are relatively easy to establish.

\begin{theorem}
\label{thm:graphic_closure}
The class of graphic Polish groups is closed under the following operations.
\begin{enumerate}[(1)]
\item If a Polish group $\Gamma$ is the internal product $\Gamma = \Gamma_1 \Gamma_2\dots \Gamma_n$ of graphic Polish groups, then $\Gamma$ is graphic. In particular, the product of finitely many graphic Polish groups is graphic. Similarly, the semidirect product of two graphic Polish groups is graphic.
\item If a Polish group $\Gamma$ is the union $\Gamma = \bigcup_{n = 1}^\infty \Gamma_n$ of countably many graphic Polish groups, then $\Gamma$ is graphic.
\item A Polish group which is a quotient of a graphic Polish group is graphic.
\item If a Polish group $\Delta$ is a closed subgroup of a graphic Polish group $\Gamma$ and has countable index in $\Gamma$, then $\Delta$ is graphic.
\end{enumerate}
\end{theorem}

\begin{proof}
All of these except (4) are essentially immediate.

(1) Given a Borel action of $\Gamma$ on $X$, look at its restriction to each $\Gamma_i$. It generates an equivalence relation that has a Borel graphing $G_i$. Then take $G = \bigcup_{i = 1}^n G_i$.

(2) Similar to (1).

(3) If $\Delta$ is a quotient of $\Gamma$ then every equivalence relation induced by a Borel action of $\Delta$ is also induced by a Borel action of $\Gamma$ (each element of $\Gamma$ acts according to its projection onto $\Delta$). Thus if $\Gamma$ is graphic then so is $\Delta$.

(4) We recall here Mackey's inducing construction which for every Borel action of $\Delta$ gives a Borel action of $\Gamma$.

Suppose $\Delta$ acts in Borel way on $X$, inducing the equivalence relation $E$. Let $T$ be a transversal for the left cosets of $\Delta$ in $\Gamma$ such that $T$ contains the identity $1\in \Gamma$. Let $Y = X\times T$. We will define a Borel action of $\Gamma$ on $Y$ as follows.

First consider the Borel action of $\Gamma$ on $T$ defined by setting $g\cdot t = s$ where $s$ is the unique element of $T$ in the coset $gt\Delta$. Next, define the cocycle $\rho\colon \Gamma\times T\to \Delta$ of this action by
\[
gt = (g\cdot t) \rho(g,t).
\]
Finally the action of $\Gamma$ on $Y$ is given by
\[
g\cdot(x,t) = (\rho (g,t)\cdot x, g\cdot t).
\]
Let $F$ be the corresponding equivalence relation.

Observe now that for $t\in T$ and $x, y\in X$,
\[
t\cdot(x,1) = (x,t)
\]
and 
\[
(x,1) F (y,1) \iff x E y.
\]
From this it follows easily that for $x,y \in X, s,t\in T$,

\[
(x,t) F (y,s) \iff x E y,
\]
and therefore $F=E\times I_T$, where $I_T = T^2$. Given now a Borel graphing $G$ of $F$ its projection to $X^2$ is a Borel graphing of $E$. Note that the fact that the projection is Borel follows from the Lusin-Novikov Theorem because $T$ is countable (and hence $T^2$ is as well).
\end{proof}

Theorem~\ref{thm:graphic_closure} above raises a few questions. First, part (4) of the theorem shows that a countable index, closed subgroup of a graphic Polish group is graphic. It seems natural to ask whether the restriction to countable index is necessary---in other words, is every closed subgroup of a graphic Polish group also graphic? We do not know the answer to this question, but we will see later that if the answer is yes, then in fact all Polish groups are graphic (see the discussion preceding Question~\ref{question:closed_subgroup}).

Second, part (1) of Theorem~\ref{thm:graphic_closure} states that any finite product of graphic Polish groups is graphic. It seems natural to ask whether this can be extended to infinite products. We will now show that this is the case, though the proof is substantially more involved.

\begin{theorem}
The class of graphic Polish groups is closed under infinite products.
\end{theorem}

\begin{proof}
Suppose that $\Gamma_0, \Gamma_1, \ldots$ is a countable sequence of graphic Polish groups and let $\Gamma = \Pi_n \Gamma_n$. Suppose we have a continuous action of $\Gamma$ on a Polish space $X$ and let $E$ denote the orbit equivalence relation of this action.

Note that given an element $g = (g_0, g_1, g_2, \ldots)$ of $\Gamma$ and a number $n$, we can naturally view $g$ as consisting of two components: namely $(g_0, g_1, \ldots, g_{n - 1})$, an element of $\Pi_{< n}\Gamma_i$, and $(g_n, g_{n + 1}, \ldots)$, an element of $\Pi_{i \geq n}\Gamma_i$. We will often need to view elements of $\Gamma$ in this way, so we now fix some notation to make this convenient. For each $n$, let $\Gamma_{< n} = \Pi_{i < n} \Gamma_i$ and $\Gamma_{\geq n} = \Pi_{i \geq n}\Gamma_i$. Note that for each $n$, we can naturally embed $\Gamma_{< n}$ into $\Gamma$ by sending $(g_0, \ldots, g_{n - 1})$ to $(g_0, \ldots, g_{n - 1}, 1, 1, \ldots)$ and we can naturally embed $\Gamma_{\geq n}$ into $\Gamma$ in a similar way. Furthermore, $\Gamma_{< n}\times \Gamma_{\geq n}$ is naturally isomorphic to $\Gamma$ via these embeddings. 

From this point on, we will identify $\Gamma_{< n}$ and $\Gamma_{\geq n}$ with their images in $\Gamma$. A basic observation that will be very important for us is that if $k_0, k_1, k_2, \ldots$ is a sequence of group elements such that for each $n$, $k_n \in \Gamma_{\geq n}$ then the infinite product $\cdots k_2k_1k_0$ is well-defined (i.e.\ the sequence of finite products $k_0, k_1k_0, k_2k_1k_0, \ldots$ converges).

We will now define the graph which witnesses graphability of $E$. To begin, fix a countable base $U_0, U_1,\ldots$ for the topology of $X$ and an oracle $a$ which knows enough about this sequence and how $\Gamma$ acts on it.\footnote{Here, ``knows enough'' essentially means ``knows whatever is needed to complete the proof'' but it is sufficent for $a$ to know the following. First, for every $i$ and $j$, $a$ knows whether the closures of $U_i$ and $U_j$ are disjoint. Second, for every $i$ and $j$, $a$ knows whether $U_i \subseteq U_j$. Third, given a point $x \in X$ and a group element $g \in \Gamma$, the list of indices $i$ such that $g\cdot x \in U_i$ should be computably enumerable in $a$ together with (codes for) $x$ and $g$. Fourth, $a$ knows countable dense subsets of each $\Gamma_{\geq n}$ and for each $k$ in one of these countable dense subsets and each $U_i$ and $U_j$, $a$ knows whether $k\cdot U_i \subseteq U_j$.} Also, for each $n$, fix a graphing $G_n$ of the action of $\Gamma_{< n}$ on $X$ (note that $\Gamma_{< n}$, being a finite product of graphic Polish groups, is also graphic). Now define a graph on $X$ as follows: two points $x, y \in X$ are adjacent in the graph if any of the following conditions holds:
\begin{itemize}
\item For some $n$, $x$ and $y$ are adjacent in $G_n$.
\item For some $n$, $x,y$ are in the same $\Gamma_{\geq n}$-orbit and that orbit has an isolated point.
\item One jump of (codes for) $x$ and $y$, together with $a$, can compute a group element $g \in \Gamma$ such that $g\cdot x = y$.
\end{itemize}
To see why the second condition above is Borel, suppose $D$ is a countable dense subset of $\Gamma_{\geq n}$ and note that the condition holds if and only if there are $g, h \in D$ and $i$ such that $g\cdot x = h\cdot y$ and for all $k \in D$, $k\cdot x \in U_i \iff k\cdot x = g\cdot x$.

Now suppose that $x$ and $y$ are in the same $\Gamma$-orbit. We will try to construct $g \in \Gamma$ such that if we define $z = g\cdot x$ then one jump of $z$ and $a$ can compute a group element taking $z$ to $x$ and a group element taking $z$ to $y$ (and hence $z$ is adjacent to both $x$ and $y$). We will show that if this fails then $x$ and $y$ are connected by a sequence of edges from the first two bullet points above.

To construct $g$, we will pick a sequence of group elements $k_0, k_1, \ldots$ such that for each $n$, $k_n$ is an element of $\Gamma_{\geq n}$. We will then define $g = \cdots k_2\cdot k_1\cdot k_0$ (recall from above that such a product is always well-defined). We will also pick a sequence of nonempty open sets $W_0 \supseteq W_1 \supseteq \ldots$ in $X$ such that for each $n$, $k_n\cdots k_0\cdot x$ is inside $W_n$. Note that since the sets $W_n$ are decreasing, $g\cdot x$ will be in the closure of each $W_n$. We will also ensure that each $W_n$ is equal to some $U_i$ (i.e.,\ some set in the countable base for the topology of $X$ that we fixed above).

The basic idea of the construction is that we will code information into our choice of the $W_n$'s. To decide what bit of information to encode on each step of the construction (i.e.\ into each choice of $W_n$) we will use a ``queue of bits'' approach similar to the proof of Theorem~\ref{thm:isomorphism}. We will explain this in more detail below.

\medskip\noindent\textit{The construction.} We will now give the details of the construction. Suppose that we have already picked $k_0, \ldots, k_n$ and $W_0, \ldots, W_n$; we will explain how to pick $k_{n + 1}$ and $W_{n + 1}$. Let $g_n = k_n \cdots k_0$. We now break into two cases.

\medskip\noindent\textit{Case 1: $g_n\cdot x$ is not the only point in its $\Gamma_{\geq n + 1}$-orbit which is inside $W_n$.} In this case, we can continue the construction and encode one more bit of information. To do so, let $(i, j)$ be the least pair of numbers such that $U_i$ and $U_j$ are subsets of $W_n$ whose closures are disjoint and which both contain some point in the $\Gamma_{\geq n + 1}$-orbit of $g_n\cdot x$. If the next bit to encode is a $0$ then set $W_{n + 1} = U_i$ and pick $k_{n + 1}$ such that $k_{n + 1}g_n\cdot x \in U_i$. Otherwise set $W_{n + 1} = U_j$ and pick $k_{n + 1}$ such that $k_{n + 1}g_n \cdot x \in U_j$.

\medskip\noindent\textit{Case 2: $g_n\cdot x$ is the unique point in its $\Gamma_{\geq n + 1}$-orbit which is inside $W_n$.} In this case, we cannot continue the construction. Instead, we will show directly that $x$ and $y$ are connected by a sequence of edges in the graph we defined above. Let $h$ and $k$ be the projections of $g_n$ onto $\Gamma_{< n + 1}$ and $\Gamma_{\geq n + 1}$, respectively. Let $g'$ be an element of $\Gamma$ such that $g'\cdot x = y$ and let $h'$ and $k'$ be the projections of $g'$ onto $\Gamma_{< n + 1}$ and $\Gamma_{\geq n + 1}$. We can now show that $x$ and $y$ are connected in the graph as follows.
\begin{enumerate}
\item First, $x$ is connected to $h\cdot x$ by a finite number of edges in the graph $G_{n + 1}$ (this follows simply from the fact that $h \in \Gamma_{< n + 1}$ and so $x$ and $h\cdot x$ are in the same $\Gamma_{< n + 1}$-orbit).
\item Next, $h\cdot x$ is connected to $k'h\cdot x$ by an edge. To see why, note that both $h\cdot x$ and $k'h\cdot x$ are in the same $\Gamma_{\geq n + 1}$-orbit as $g_n\cdot x = kh\cdot x$ and by assumption, $g_n\cdot x$ is the unique point in its $\Gamma_{\geq n + 1}$-orbit which is inside the open set $W_n$.
\item Finally, $k'h\cdot x$ is connected to $h'h^{-1}k'h\cdot x$ by a finite number of edges in the graph $G_{n + 1}$ (this just follows from the fact that $h'h^{-1} \in \Gamma_{< n + 1}$). Since $\Gamma_{< n + 1}$ and $\Gamma_{\geq n + 1}$ commute with each other, we have
\[
  h'h^{-1}k'h = k'h'h^{-1}h = k'h' = g'
\]
and thus $h'h^{-1}k'h\cdot x = g'\cdot x = y$.
\end{enumerate}

\medskip\noindent\textit{Verification.} We will now assume that we never encounter case 2 in the construction above, in which case after infinitely many stages the construction builds an element $g = \cdots k_2k_1k_0$ of $\Gamma$. Let $z = g\cdot x$. We will first explain why one jump of $z$ and $a$ can compute the information encoded during the construction of $g$. We will then explain what information should be encoded.

To see why the information can be decoded, we will first explain how to compute the sequence of open sets $W_0, W_1, W_2, \ldots$. Suppose that we have successfully determined $W_n$ and we want to determine $W_{n + 1}$. First, find the least pair $(i, j)$ such that $U_i$ and $U_j$ are subsets of $W_n$ with disjoint closures and which both contain an element of the $\Gamma_{\geq n + 1}$-orbit of $z$. Next, check whether $z$ is in the closure of $U_i$ or $U_j$. The one whose closure contains $z$ is $W_{n + 1}$. Now note that when we determine whether $W_{n + 1} = U_i$ or $W_{n + 1} = U_j$, we also find out the bit that was encoded on step $n$ of the construction---in particular, if $W_{n + 1} = U_i$ then the bit encoded was a $0$ and otherwise it was a $1$.

Now we will explain what information should be encoded. Fix a group element $g'$ such that $g'\cdot x = y$. It is enough to encode $g'$ along with the sequence $k_0, k_1, k_2, \ldots$ chosen during the construction (because if we know this sequence then we can compute $g$ and thus we know group elements taking $x$ to $z$ and $x$ to $y$). As in the proof of Theorem~\ref{thm:isomorphism}, there is an apparent circularity: the choice of the $k_n$'s depends on what information we need to encode, which in turn depends on the choice of the $k_n$. To get around this problem, we will use a variation on the ``queue of bits'' idea from the proof of Theorem~\ref{thm:isomorphism}. Imagine forming an infinite table. Fill in the first row of this table with bits encoding $g'$. Fill in the rest of the rows during the course of the construction, on step $n$ adding bits encoding $k_n$ to row $n + 1$. Now pick some path which visits every cell of this table and on step $n$ of the construction, encode the bit in the $n^\text{th}$ cell visited by this path. As long as the path doesn't visit row $n$ before step $n$, there is no circularity.
\end{proof}

\subsection{Connected groups}

We will now show that all connected Polish groups are graphic. We will prove this in two steps. First, we will formulate a weakening of the Borel witness coding property, which we call the Borel coding property, and prove that for equivalence relations arising from Polish group actions, this weaker coding property is still enough to imply Borel graphability. Second, we will show that the orbit equivalence relation of any action of a connected group has this coding property.

\begin{definition}
\label{def:borelcoding}
An analytic equivalence relation $E$ on $X$ has the \term{Borel coding property} if there is a Borel function $f\colon X\to \Cantor$ such that for all $x \in X$ and $a \in \Cantor$, there is some $y$ which is $E$-equivalent to $x$ such that $f(y)$ computes $a$.
\end{definition}

Note that this is almost the same as the Borel witness coding property, except that $y$ is not required to encode a witness to the equivalence $xEy$.

Our first goal is to show that if $E$ is induced by a Polish group action and has the Borel coding property then $E$ is Borel graphable. To prove this, we will use the following technical lemma.

\begin{lemma}
Let $\Gamma$ be a Polish group acting continuously on a Polish space $X$ and suppose that every orbit is uncountable. Then there is a Borel map $ \gamma \colon X \times \Cantor\to \Gamma$ such that for every $x\in X$ and $a \not= b \in \Cantor$, $\gamma (x,a)\cdot x \not= \gamma (x,b)\cdot x$.
\end{lemma}

Before giving the proof of this lemma, we will briefly explain its meaning. Suppose that $\Gamma$ is a Polish group acting continuously on a Polish space $X$ and that every orbit is uncountable. Now fix some $x \in X$. Since the orbit of $x$ is uncountable, it follows from the perfect set theorem that it contains a perfect set, $A_x$. Moreover, it follows from the proof of the perfect set theorem that it is actually possible to find a perfect set $B_x$ in $\Gamma$ which generates this perfect set---in particular, such that the map $g \mapsto g\cdot x$ is injective on $B_x$ and its image is $A_x$. The lemma above essentially tells us that we can pick the sets $B_x$ in such a way that the map $x \mapsto B_x$ is Borel.

\begin{proof}
Let $K(\Gamma)$ be the Polish space of compact subsets of $\Gamma$ (see Section 4.F of \cite{kechris1995classical}). Define $R\subseteq X \times K(\Gamma)$ by
\[
R(x,K) \iff \forall g,h (g,h \in K \ \& \ g\not= h \implies g\cdot x \not= h\cdot x).
\]
We first note that $R$ is Borel, in fact $G_\delta$. To see this, let $\{U_n\}$ be an open basis for $G$ and for each $n, m$ such that $\bar{U_m}\cap \bar{U_n} = \0$, define 
\[
    S_{m, n}(x, K) \iff \exists g, h~(g \in K \cap \bar{U_m} \ \&\ h \in K \cap \bar{U_n} \ \&\ g\cdot x = h\cdot x). 
\]
Note that each $S_{m, n}$ is closed and that
\[
    \neg R(x, K) \iff \exists m, n~(\bar{U_m} \cap \bar{U_n} = \0 \ \&\ S_{m, n}(x, K))
\]
and hence $\neg R$ is $F_\sigma$.

Next we show that for each $x$, $R_x$ is comeager. To see this let
\[
L_x = \{(g,h)\colon g\cdot x = h\cdot x\}.
\]
Note that for each $g$, $(L_x)_g = g\Gamma_x$ where $\Gamma_x$ is the stabilizer of $x$. Since the orbit of $x$ is uncountable, $\Gamma_x$ has uncountable index in $\Gamma$ and is thus meager and so $(L_x)_g$ is meager as well. Therefore $L_x$ is meager by the Kuratowski-Ulam Theorem (see 8.41 in~\cite{kechris1995classical}). 
Since 
\[
R_x = \{ K\in K(\Gamma)\colon \forall g,h (g,h \in K \ \&  \ g\not= h \implies g\cdot x \not= h\cdot x)
\},
\]
it follows from the Kuratowski, Mycielski Theorem (see 19.1 in \cite{kechris1995classical})
that $R_x$ is comeager.

Next, note that since $R_x$ is comeager, so is its intersection with the perfect sets in $K(\Gamma)$, which form a dense $G_\delta$ set in $K(\Gamma)$ (see 8.8 in \cite{kechris1995classical}). So by the uniformization theorem for Borel sets with non-meager sections, there is a Borel function $k\colon X \mapsto K(\Gamma)$ such that $\forall x [R(x, k(x)) \ \textrm{and} \ k(x) 
 \ \textrm{is perfect}]$. 
 
Finally, by standard arguments, there is a Borel function $f\colon K(\Gamma)\times \Cantor\to \Gamma$ such that if $K$ is perfect, then $f(K,a)\in K$ and $a\not=b\implies f(K,a) \not= f(K,b)$. We can now finish by setting $\gamma(x,a) = f(k(x), a)$.
\end{proof}

Our use of the above lemma will be encapsulated in the following corollary. Since the statement of the corollary may look somewhat obscure, we will explain how one can think about it informally. It essentially states that given any $x \in X$ and $a \in \Cantor$, we can encode $a$ into some element $y$ of the orbit of $x$ such that $a$ can be decoded using $x$ and $y$ together. The function $e$ in the statement of the corollary accomplishes the encoding (i.e.\ it tells us how to find $y$), the function $d$ accomplishes the decoding and the function $g$ provides a witness that $x$ and $y$ are in the same orbit.

\begin{corollary}
\label{cor:polish_group_pair_coding}
Let $\Gamma$ be a Polish group acting continuously on a Polish space $X$ and suppose that every orbit is uncountable. Then there are Borel functions $e\colon X\times X \to \Cantor$, $d\colon X\times X \to \Cantor$ and $g \colon X\times X \to \Gamma$ such that for all $x \in X$ and $a \in \Cantor$, if $y = e(x, a)$ then 
\begin{enumerate}
    \item $d(x, y) = a$
    \item and $y$ is in the same orbit as $x$, as witnessed by $g(x, y)$---i.e.\ $g(x, y) \cdot x = y$.
\end{enumerate}
\end{corollary}

\begin{proof}
Let $\gamma \colon X \times \Cantor \to \Gamma$ be the function given by the lemma above. Define $e(x, a) = \gamma(x, a)\cdot x$. Note that for a fixed $x$, the function $a\mapsto e(x, a)$ is injective. Thus by the Lusin-Novikov Theorem, there is a Borel function $d\colon X\times X \to \Cantor$ such that for all $x \in X$ and $a \in \Cantor$, $d(x, e(x, a)) = a$. Now define $g(x, y) = \gamma(x, d(x, y))$. It is straightforward to check that these functions have the desired properties.
\end{proof}

\begin{theorem}
\label{thm:coding_to_graphable}
Let $\Gamma$ be a Polish group acting continuously on a Polish space $X$ and let $E$ be the associated orbit equivalence relation. If $E$ has the Borel coding property, then it is Borel graphable (with diameter $4$).
\end{theorem}

\begin{proof}
Note that if $E$ has the Borel coding property then all orbits must be uncountable. Thus Corollary~\ref{cor:polish_group_pair_coding} applies so we have functions $e, d,$ and $g$ with the properties described in that corollary. Let $f\colon X \times \Cantor \to X$ be a function witnessing that $E$ has the Borel coding property.

For convenience, we fix Borel bijections of $X$ and $\Gamma$ with $\Cantor$ so that we may assume $X = \Gamma = \Cantor$ (without affecting the fact that the group operations of $\Gamma$ and the action of $\Gamma$ on $X$ are both Borel). Given $x,y\in X$, let $A_{x,y}$ be the smallest subset of $\Cantor$ containing $x,y$ and closed under the group operations, $d,g, f$, finite joins (in the sense of computability  theory) and Turing reductions. Define the graph $G$ by 
\[
x G y \iff \exists g \in A_{x,y}~ (g\cdot x = y \text{ or } g\cdot y = x).
\]
Clearly $G$ is Borel.

Now suppose $x$ and $y$ are in the same orbit. We will show they are connected in the graph described above by a path of length $4$. Using the Borel coding property, pick an element $w$ of the orbit of $x$ such that $f(w)$ computes $x$ and $y$. Let $a_1$ and $a_2$ be elements of $\Gamma$ such that $a_1\cdot w = x$ and $a_2\cdot w = y$. Let $z_1 = e(x, a_1)$ and $z_2 = e(y, a_2)$.

We claim that $x, z_1, w, z_2, y$ forms a path in the graph defined above. It is enough to check that there is an edge between $x$ and $z_1$ and between $z_1$ and $w$; the cases of $w, z_2$ and $z_2, y$ are identical. To see why there is an edge between $x$ and $z_1$, simply note that $g(x, z_1)$ is a group element taking $x$ to $z_1$. To see why there is an edge between $z_1$ and $w$, we can reason as follows.
\begin{enumerate}
\item $f(w)$ computes $x$ by our choice of $w$.
\item $d(x, z_1)$ computes $a_1$.
\item $g(x, z_1)a_1$ is a group element taking $w$ to $z_1$. \qedhere
\end{enumerate}
\end{proof}

We are now ready to apply Theorem \ref{thm:coding_to_graphable} to prove
that all connected Polish groups are graphic.

\begin{theorem}
\label{thm:connected}
Suppose that $\Gamma$ is a connected Polish group acting continuously on a Polish space $X$ and let $E$ be the associated orbit equivalence relation. If there are no fixed points for this action, then $E$ has the Borel coding property.
\end{theorem}

\begin{proof}
The main idea is to code information into distances. Fix a countable dense set $\{w_i\}_{i \in \N}$ of elements of $X$ and a metric $d$ on $X$. Define $f \colon X \to \Cantor$ to be the join of codes of $d(x, w_i)$ over all $i \in \N$. Now fix $x \in X$ and $a \in \Cantor$; we will show that for some $z$ in the orbit of $x$, $f(z)$ computes $a$.

We claim that there is some $y$ in the orbit of $x$ and some $i$ such that $d(x, w_i) \neq d(y, w_i)$. If not then we must have that $x$ is a fixed point of $\Gamma$, contradicting our assumption.

Since $g\mapsto d(g\cdot x, w_i)$ is continuous, its range is a non-trivial interval in $\R$. So there is $g\in \Gamma$ with the code of $r= d(g\cdot x, w_i)$ computing $a$. Let $z=g\cdot x.$
\end{proof}

\begin{corollary}
Every connected Polish group is graphic.
\end{corollary}

\begin{proof}
Consider a continuous action of a connected Polish group $\Gamma$ on a Polish space $X$. The set of fixed points of this action is closed, so it is enough to consider the action on the open set $Y\subseteq X$ of non-fixed points. The corollary then follows from Theorems~\ref{thm:coding_to_graphable} and~\ref{thm:connected}.
\end{proof}

\begin{remark}
Note that if a connected Polish group $\Gamma$ acts continuously on a Polish space $X$ and $x \in X$ is not a fixed point, then the orbit $\Gamma\cdot x$ is uncountable. This is because the stabilizer $\Gamma_x$ is not open, so it has uncountable index.     
\end{remark}

\begin{remark}
The above proof does not show that every orbit equivalence relation of a connected Polish group action has a Borel graphing of diameter $2$. However, it can be shown by a similar proof that if a Polish group is path connected, then all of its orbit equivalence relations admit Borel graphings of diameter $2$.
\end{remark}

Recall that we began this section by asking whether every Polish group is graphic. The fact that all connected Polish groups are graphic has an interesting implication for this question. Hartman and Mycielski~\cite{hartman1958} proved that every Polish group is a closed subgroup of a connected Polish group. Thus an affirmative answer to the following question would imply that all Polish groups are graphic.

\begin{question}
\label{question:closed_subgroup}
Are the closed subgroups of a graphic Polish group  graphic?
\end{question}

\subsection{Non-archimedean groups}
\label{sec:nonarchimedean}

A Polish group is called non-archimedean if it is isomorphic (as a topological group) to a closed subgroup of $S_\infty$. Equivalently, a Polish group is non-archimedean if it is isomorphic to the automorphism group $\Aut(M)$ of some countable ultrahomogeneous structure $M$ in a countable relational language. Recall that a structure $M$ is \term{ultrahomogeneous} if every isomorphism between finite subsets of $M$ can be extended to an automorphism of $M$ (for more about non-archimedean Polish groups, see the survey~\cite{kechris2013dynamics}).

One motivation for investigating graphability for non-archimedean groups comes from a result of Ding. A Polish group $\Gamma$ is \term{projectively universal} if every Polish group $\Delta$ is a quotient of $\Gamma$ by a closed subgroup. If one could show that some projectively universal Polish group is graphic, then it would follow by Theorem~\ref{thm:graphic_closure} (3) that every Polish group is graphic. Ding~\cite{ding2012} showed that projectively universal Polish groups exist and, moreover, there is a projectively universal group which is totally disconnected. Thus, to prove that all Polish groups are graphic, it would be enough to prove that all totally disconnected Polish groups are graphic. Since the non-archimedean groups are standard examples of totally disconnected Polish groups, it is reasonable to first try to prove that they are all graphic.

We are not quite able to prove that all non-archimedean groups are graphic, but our next result states that a fairly large subset of them are. In particular, if $\Gamma$ is the automorphism group of a countable ultrahomogeneous structure in a \emph{finite} relational language, then $\Gamma$ is graphic.

To prove this, we will use a result about universal actions for non-archimedean groups due to Becker and Kechris~\cite{becker1996descriptive} (which generalizes the result about $S_\infty$ mentioned at the beginning of this section). Fix a non-archimedean Polish group $\Gamma$ and let $\lang$ be a countable relational language and $A$ an ultrahomogeneous $\lang$-structure with domain $\N$ such that $\Gamma$ is isomorphic to $\Aut(A)$. Now fix a relational language $\lang'$ with relations of arbitrarily large arity and consider the equivalence relation on $\Mod(\lang')$ given by setting $\lang'$-structures $M$ and $N$ equivalent if they are isomorphic via a function $f\colon \N \to \N$ which is an automorphism of $A$. This equivalence relation is naturally induced by an action of $\Gamma$ and is universal for Polish group actions of $\Gamma$. In particular, any orbit equivalence relation of a Polish group action of $\Gamma$ has an invariant Borel reduction to this equivalence relation.

\begin{theorem}
Suppose $\Gamma$ is the automorphism group of a countable ultrahomogeneous structure in a finite relational language. Then $\Gamma$ is graphic.
\end{theorem}

\begin{proof}
Fix a finite relational language $\lang$ and an ultrahomogeneous $\lang$-structure $A$ with domain $\N$ such that $\Gamma = \Aut(A)$. Also fix a countable relational language $\lang'$ with relations of arbitrarily large arity and let $E$ be the equivalence relation on $\Mod(\lang')$ discussed above. It is enough to show that $E$ is Borel graphable. By Proposition~\ref{prop:countable_orbit1}, we can restrict our attention to the uncountable equivalence classes of $E$.

In fact, we will show that $E$ restricted to the uncountable equivalence classes has the Borel witness coding property. In particular, we will show that for any $\lang'$-structure $M$ with domain $\N$ with uncountable $E$-equivalence class and any real $r \in \Cantor$, there is some $\lang'$-structure $N$ with domain $\N$ and some $f \in \Aut(A)$ which is an isomorphism from $N$ to $M$ and such that a Turing jump of $N\oplus A\oplus \lang'$ computes $f$ and $r$.

Our proof is similar to our proof of Theorem~\ref{thm:isomorphism}, but with a few extra complications along the way. As in that proof, we will use the mental image of constructing $N$ by placing elements of $M$ on natural numbers. Recall that, formally, this corresponds to building a bijection $f \colon \N \to M$ and defining $N$ by setting, for each $\lang'$-relation $R$,
\[
    R^N(n_0, \ldots, n_k) \iff R^M(f(n_0), \ldots, f(n_k)).
\]
Also recall that placing an element $a$ of $M$ on $n$ corresponds to setting $f(n) = a$. 

Also as in the proof of Theorem~\ref{thm:isomorphism}, we will build $N$ in stages. On each stage, we will place some more elements of $M$ in order to code one more bit of information and then place one additional element to ensure that every element of $M$ is placed somewhere. Unlike in the proof of Theorem~\ref{thm:isomorphism}, however, we need to make sure that the function $f$ we build is an automorphism of $A$. To do so, it is enough to ensure that at each stage, the finite fragment of $f$ we have constructed so far can be extended to an automorphism of $A$.

Some terminology will be useful when discussing this. Say that a placement of elements of $M$ on natural numbers is \term{safe} if the corresponding function can be extended to an automorphism of $A$. Suppose we have a safe placement of elements $a_0, \ldots, a_n$ on the numbers $0, \ldots, n$ and we want to place an element $b$ on a number $k$. The fact that $A$ is ultrahomogeneous implies that this new placement is still safe as long as $b$ has the same quantifier-free $\lang$-type over $a_0,\ldots,a_n$ that $k$ has over $0, \ldots, n$. 

\medskip\noindent\textit{The construction.} We will now give the details of the construction. Suppose we have just finished stage $n$ and have placed elements $a_0,\ldots,a_k$ of $M$ on numbers $0,\ldots,k$. Also assume that this placement is safe. We will explain what to do on stage $n + 1$. We first claim that for some $l$, there are length $l$ tuples $\bar{b}$ and $\bar{c}$ in $M$ consisting of unplaced elements which can be safely placed on the next available numbers (i.e.\ on the numbers $k + 1,\ldots, k + l$) and which have different quantifier-free $\lang'$-types over $\{a_0,\ldots,a_k\}$. If not, then it is not hard to show that the $G$-orbit of $M$ is countable (essentially by following the proof of Lemma~\ref{lemma:trivial_isomorphism1}).

Now let $\phi$ be the least quantifier free $\lang'$-formula (with the least number of free variables) for which there are tuples $\bar{b}$ and $\bar{c}$ as described above for which $M \models \phi(\bar{a}, \bar{b})$ and $M \notmodels \phi(\bar{a}, \bar{c})$. As in the proof of Theorem~\ref{thm:isomorphism}, we now code one bit of information by our choice of which of these tuples to place next: if the next bit of information to code is a $0$ then we place the elements of $\bar{b}$ next; otherwise we place the elements of $\bar{c}$.

To finish the current stage, we need to place one more element of $M$ in order to make sure that all elements are eventually placed somewhere. This is the only part of the proof that requires more care than the proof of Theorem~\ref{thm:isomorphism}. The basic issue is as follows. Let $b$ be the least element of $M$ not yet placed somewhere and recall that $k + l + 1$ is the first available number which does not have something placed on it. Unlike in the previous proof, we cannot simply place $b$ on $k + l + 1$ since the resulting placement may not be safe. On the other hand, if we place $b$ somewhere else, then this may cause problems for the decoding process (since on each step of the decoding, we need to know which numbers have already had elements of $M$ placed on them by the beginning of the corresponding step of the construction). We will get around this problem by using the fact that $\lang$ is a finite language.

Note that since $\lang$ is finite, there are only finitely many possible quantifier-free $\lang$-types over the elements of $M$ that we have placed so far. Suppose there are $p$ such possible types and order them lexicographically. Now for each $i \leq p$, let $m_i$ denote the least natural number which realizes the $i^\text{th}$ of these types over the numbers which have had elements placed on them already (i.e. over the numbers $0, \ldots, k + l$). The point is that for each $c$ in $M$ which has not yet been placed, there is some $m_i$ on which we can safely place $c$, given the placements we have already made. Let $m = \max_{i \leq p}m_i$. Place elements of $M$ on all numbers up to $m$, making sure to place $b$ on one of these elements (note that we can do this by first safely placing $b$ on the $m_i$ whose type matches that of $b$ and then choosing an arbitrary safe placement for all the remaining numbers up to $m$).

\medskip\noindent\textit{Verification.} The key point is that doing things in this way ensures that during the decoding process, we can figure out which numbers had elements placed on them on each stage. We will now explain this point in more detail; the rest of the proof (in particular, how to decide what information to encode on each stage) is almost identical to the proof of Theorem~\ref{thm:isomorphism}.

Suppose that we are trying to use $N$ to decode what information was encoded on each stage and we know that by the end of stage $n$, elements were placed on the numbers $0, \ldots, k$. We first look for the least $l$ such that there is some tuple which has the same quantifier-free $\lang$-type over $\{0, \ldots, k\}$ as $(k + 1, \ldots, k + l)$, but a different $\lang'$-type. This tells us the bit encoded on stage $n + 1$ as well as how many elements were placed in the first part of the stage. Then, for each quantifier-free $\lang$-type over $\{0, \ldots, k + l\}$, we look for the least number which realizes this type and we let $m$ be the maximum of these numbers. This tells us which numbers had elements placed on them in the second part of stage $n + 1$ (in particular, it tells us that elements were placed on all numbers up to $m$).
\end{proof}


\section{Structural properties of Borel graphability}
\label{sec:closure_properties}

In this section, we will address various structural properties
of the class of Borel graphable equivalence relations.

\subsection{Effectivity}

We will establish here that Borel graphability is not effective, in the sense that the existence of a Borel graphing for a $\Sigma^1_1$ equivalence relation does not imply the existence of a $\Delta^1_1$-graphing (a fact first proved by the first author in~\cite{arant2019effective}).  Note that assuming there is a non-constructible real, this follows from Propositions~\ref{prop:not_delta11} and~\ref{prop:graphable}: $\Fw$ is $\Sigma^1_1$ and not $\Delta^1_1$-graphable, but is Borel graphable if there is a non-constructible real.
However, using the results above it is also easy to find
an example establishing the non-effectivity of Borel
graphability that works in ZFC.

\begin{theorem}
\label{thm:noneffective}
There is a $\Sigma^1_1$ equivalence relation which is Borel graphable but not $\Delta^1_1$-graphable.
\end{theorem}

\begin{proof}
Let $X_\text{low}$ be the set of all reals $x\in \Cantor$
with $\omega_1^x=\omega_1^{\text{CK}}$.  This set is both
$\Sigma^1_1$ and $\Delta^1_1(\mathcal{O})$ (see
\cite{chong2015}, Section 2.4).  Thus, the equivalence
relation $E$ on $\Cantor$ given by
\[
xEy \iff x, y \in X_{\text{low}} \ \text{or} \ x=y
\]
is both $\Sigma^1_1$ and Borel.  In particular,
$E$ is Borel graphable.  However, $E$ is not
$\Delta^1_1$-graphable.  To see this, first note
that $E$ has
$X_\text{low}$ as an equivalence class.  The proof
of Proposition \ref{prop:not_delta11} actually establishes
that if an equivalence relation $E$ has $X_{\text{low}}$
as an equivalence class, then $E$ is not 
$\Delta^1_1$-graphable.  
\end{proof}

\subsection{Closure properties for Borel graphability}

We have previously seen in Corollary \ref{cor:borelreducibility} and Proposition \ref{prop:invariantreducability}
that Borel graphability is not closed under Borel reducibility,
but is closed under \textit{invariant} Borel reducibility.  
The following closure properties are also easily established.

\begin{proposition}
The class of Borel graphable equivalence relations is closed under the following operations.
    \begin{itemize}
        \item[(1)] Finite products.
        \item[(2)] Countable direct sums.
    \end{itemize}
\end{proposition}

\begin{proof}
    For (1), if $E_i$ is Borel graphed by $G_i$ for $i<n$, 
    just define
    \[
    (x_0, \dots, x_{n-1}) G (y_0, \dots, y_{n-1}) \iff 
    (\forall i < n)~(x_iG_i y_i \text{ or } x_i = y_i)
    \]

    For (2), just take the countable direct sum of the Borel graphings 
    of the equivalence relations $E_i$.
\end{proof}

Item (1) from the previous proposition naturally raises the
question of whether Borel graphability is closed under taking countably infinite products. We will show that this is true, but unlike the case for finite products, the proof is non-trivial.

\begin{theorem}
\label{thm:infproducts}
Borel graphability is closed under infinite products.
\end{theorem}

\begin{proof}
Without loss of generality, suppose that each $E_i$
is an equivalence relation on $\Cantor$ and set $E = \prod_i E_i$, which is an equivalence relation on $X = (\Cantor)^\N$.

For each $i\in \N$, let $C_i$ be an unfolding of $E_i$ (see Definition~\ref{def:unfolding}) and let $G_i$ be a Borel graphing of $E_i$. Define a function $h:X\times X \rightarrow \Cantor$ by
\[
h(x, y)(i) = \begin{cases} 1 & \text{if $x_iG_iy_i$} \\
                    0 & \text{otherwise}, \end{cases}
\]
which is clearly Borel.

Our Borel graphing $G$ of $E$ will have two types of edges
between distinct $x, y\in X$:

\begin{enumerate}
    \item an edge between $x$ and $y$ 
    if $(\forall i)[x_i=y_i \vee x_i G_i y_i]$; and
    
    \item an edge between $x$ and $y$ if 
    a jump of $x\oplus y \oplus h(x, y)$ computes a witness
    to $x_iE_iy_i$ for every $i\in \N$. Recall that by ``a witness to $x_i E_i y_i$, we mean a real $w \in \Cantor$ such that $(x_i, y_i, w) \in C_i$.
\end{enumerate}

Note that if there is a uniform bound for the diameters of the $G_i$,
we can obtain a graphing for $E$ just using edges of type (1).  
We will use the notation $\text{dist}_i(x, y)$ 
to denote the distance between $x$ and $y$ in the graph $G_i$.

We now show that $G$ is a graphing of $E$.  
Let $x, y\in X$ be distinct with $x E y$.  
By using the edges of type (1), we can find a path from $x$ to some
$z$ such that, for all $i\in \N$, 
either $z_i=y_i$ or $\text{dist}_i(z_i, y_i)\geq 3$.  
Let $I=\{i : \text{dist}_i(z_i, y_i)\geq 3\}$.  
We may assume $I$ is infinite, since otherwise
we can obviously find a path from $z$ to $y$ using edges of type (1).

We will now describe how to construct 
$w\in X$ so that both
$(z\oplus w \oplus h(z, w))'$ and  $(y\oplus w \oplus h(y, w))'$ compute 
an arbitrary piece of information.
Later, we will carefully explain what information should be encoded to ensure that there is an edge between $z$ and $w$ and between $w$ and $y$.  
For $i\notin I$, set $w_i=z_i=y_i$ 
(the information will all be encoded in the $w_i$ for which $i\in I$).
For each $i\in I$, fix a minimal length path in $G_i$ from $z_i$ to $y_i$, 
which necessarily has length $\geq 3$.
Note that this path has a unique element $a_i$ 
which is adjacent to $z_i$ but not $y_i$ and
a unique element $b_i$ which is adjacent to $y_i$ but not $z_i$.  
To encode a $0$, set $w_i=a_i$, and to encode a $1$ set $w_i=b_i$.  
We can easily decode the construction using
$z\oplus w \oplus h(z, w)$; e.g., to check if a $1$ 
is being coded in the $i^\text{th}$ coordinate,
check that $z_i\neq w_i$ and that $z_i$ is not $G_i$-adjacent to $w_i$.
Similarly, $y\oplus w \oplus h(y, w)$ can also decode the information.

Now, we have to describe what information 
we want to encode in the construction so
that there are edges of type (2) from $z$ to $w$ 
and from $y$ to $w$.
This means we want to encode witnesses to $z_iE_iw_i$ and $y_iE_i w_i$
for all $i\in \N$. 
Obviously, which witnesses need to be encoded
is determined as we carry out the construction of $w$.  
To avoid having to constantly update what information we need to encode,
we can use the following trick.  Before we begin the construction of $w$,
for each $i\in \N$ fix a minimal length $G_i$-path 
$(x^i_0, \dots,x^i_{n_i-1})$ with $x^i_0=z_i$ and $x^i_{n_i-1}=y_i$.
Fix witnesses $c^i_{jk}$ to $x^i_jE_ix^i_k$ for all $0\leq j, k< n_i$,
and let $c$ be the join of all the $c^i_{jk}$.  Then,
we code the bits of $c$ into our construction of $w$.  It follows
that $zGw$ and $yGw$, as desired.
\end{proof}

The next result establishes a criterion for
the Borel graphability for infinite products that does
\textit{not} require that the factor equivalence relations be
Borel graphable. 

\begin{theorem}
    If $E_0, E_1, \dots$ is an infinite sequence of 
    analytic equivalence relations with all equivalence
    classes of size $>1$, then the infinite product
    of these equivalence relations is Borel graphable.
\end{theorem}

\begin{proof}
    The proof is similar to that of Theorem \ref{thm:infproducts}.
    We may assume each $E_i$ is an equivalence relation on $\Cantor$. 
    Let $E=\prod E_i$ be the product equivalence relation on $X=(2^\N)^\N$ and for each $i$, let $C_i$ be an unfolding of $E_i$.
    The Borel graphing $G$ of $E$ is defined so that there
    is an edge between distinct  $x, y\in X$ if 
    a jump of $x\oplus y$ computes witnesses to $x_iE_iy_i$ 
    for every $i\in \N$ (where by a ``witness to $x_i E_i y_i$,'' we mean a real $w \in \Cantor$ such that $(x_i, y_i, w) \in C_i$).

    Let $x, y\in X$ be distinct with $xEy$. We will describe
    how to construct a new element $z\in X$ such that
    $(x\oplus z)'$ and $(y\oplus z)'$ both compute an arbitrary
    piece of information. By our hypothesis, for every $i\in \N$
    we may fix
    $a_i\in [x_i]_{E_i}$ with $a_i\neq x_i$ and $b_i\in [y_i]_{E_i}$
    with $b_i\neq y_i$.
    The $i^\text{th}$ bit will be encoded by both $z_{2i}$ and $z_{2i+1}$.
    If the $i^\text{th}$ bit is a $0$, set $z_{2i}=x_{2i}$ and $z_{2i+1}=y_{2i+1}$.
    If the $i^\text{th}$ bit is a $1$, set $z_{2i}=a_{2i}$ and $z_{2i+1}=b_{2i+1}$.
    It is easy to see that both $x\oplus z$ and 
    $y\oplus z$ can decode the construction; e.g.,
    using $x\oplus z$, we can check if the $i^\text{th}$ bit
    is a $0$ or a $1$ by checking whether $x_{2i}=z_{2i}$ or not.

    The information that we need to encode consists of witnesses
    to $x_iE_iz_i$ and $y_iE_iz_i$ for all $i\in\N$.  Thus, we can
    just encode the join of witnesses to $x_iE_i a_i$ (for $i$ even), $y_i E_i b_i$ (for $i$ odd) and $x_iE_i y_i$
    (for all $i\in \N$). It follows that $xGz$ and $zGy$.
\end{proof}

We end this subsection by discussing some closure properties
which fail for Borel graphable equivalence relations.

\begin{proposition}
    The class of Borel graphable equivalence relations
    is not closed under the following operations:
    \begin{enumerate}[(1)]
    \item uncountable direct sums;
    \item taking subequivalence relations;
    \item taking superequivalence relations.
    \end{enumerate}
\end{proposition}

\begin{proof}
    For (1), just note that the non-Borel graphable 
    equivalence relation in the proof
    of Proposition \ref{prop:firstnonexample} 
    is a direct sum of uncountably many
    equivalence relations on the finite space $\{0, 1\}$.

    For (2), start with a non-Borel graphable equivalence relation $E$
    on $X$.  On the space  $Y=X\cup (X\times\Cantor)$, define an
    equivalence relation $E'$ so that its classes are exactly of
    the form $C\subseteq X$ and $C\times \Cantor\subseteq X\times \Cantor$,
    where $C$ is an $E$-class.  Clearly, $E'$ is not Borel graphable.
    Now, define $E''$ on $Y$ so that its equivalence classes are
    exactly of the form $C\cup (C\times \Cantor)$, where $C$
    is an $E$-class.  $E''$ is Borel graphable by placing
    edges between $x$ and $(x, a)$ for every $x\in X$ and $a\in \Cantor$,
    and placing an edge from $(x, a)$ to $(y, a)$ whenever $a$
    is a witness to $xEy$.  Then, $E'\subseteq E''$, where $E'$
    is not Borel graphable but $E''$ is Borel graphable.  

    (3) is simple, since equality is a Borel graphable subequivalence
    relation of every equivalence relation.  
\end{proof}

\subsection{Borel coding vs. Borel graphability}
\label{sec:borelcoding}

Several times throughout this paper, we have proved that an equivalence relation $E$ is Borel graphable by showing, in essence, that arbitrary information can be encoded into elements of every (or almost every) equivalence class of $E$. Examples can be found in our proofs of Borel graphability for $\Fw$, the isomorphism relation on $\Mod(\lang)$, actions of infinite products of graphic Polish groups, actions of connected Polish groups and infinite products of Borel graphable equivalence relations. Based on these examples, it is hard to escape the feeling that, at least in many cases, Borel graphability is equivalent to this sort of ability to code information.

We have already seen two cases in which we could formally state and prove an implication from some sort of coding property to graphability. First, in Section~\ref{sec:basics} we defined the Borel witness coding property and proved that it implies Borel graphability. Second, in Section~\ref{sec:graphic} we defined a somewhat weaker property, the Borel coding property, and proved that for equivalence relations arising from Polish group actions, it implies Borel graphability as well. However, we have so far not addressed the other direction: to what extent does Borel graphability require some coding property to hold?

There is one trivial barrier to proving such an implication. Namely, if one can encode arbitrary information into the elements of each equivalence class of $E$, then presumably the equivalence classes of $E$ must all be uncountable. However, certainly not all Borel graphable equivalence relations have all equivalence classes uncountable. Thus when studying to what extent graphability implies coding, we will only work with equivalence relations whose equivalence classes are all uncountable.

With this caveat, we are able to prove that Borel graphability implies a certain weak coding property. To state our result, we need to introduce some definitions. First, we will introduce a weakening of the Borel coding property.

\begin{definition}
An equivalence relation $E$ on a Polish space $X$ has the \term{Borel pair-coding property} if there is a Borel function $f\colon X\times X \to \Cantor$ such that for all $x \in X$ and $a \in \Cantor$, there is some $y$ which is $E$-equivalent to $x$ such that $f(x, y)$ computes $a$.
\end{definition}

Note that we have now introduced three different sorts of coding properties; we will now list them from strongest to weakest. In each case, an analytic equivalence relation $E$ on $X$ has the coding property if there is a Borel function $f$ such that for all $x\in X$ and $a \in \Cantor$, there is some $y$ which is $E$-equivalent to $x$ such that:
\begin{itemize}
\item \textbf{Borel witness coding:} $f(y)$ computes both $a$ and a witness to $x E y$.
\item \textbf{Borel coding:} $f(y)$ computes $a$.
\item \textbf{Borel pair-coding:} $f(x, y)$ computes $a$.
\end{itemize}

Recall that a function $f$ is Borel if and only if it is hyperarithmetic relative to some oracle. This gives us a natural way to weaken each of the above coding properties: we can replace the function $f$ with a hyperarithmetic reduction which is \emph{allowed to depend on $y$} (where by ``hyperarithmetic'' we mean hyperarithmetic relative to a fixed oracle which does not depend on $y$). We will refer to these weakened versions as the \term{weak Borel witness coding property}, the \term{weak Borel coding property} and the \term{weak Borel pair-coding property}. For example, the weak Borel pair-coding property can be defined formally as follows.

\begin{definition}
An analytic equivalence relation $E$ on $\Cantor$ has the \term{weak Borel pair-coding property} if there is an oracle $a \in \Cantor$ such that for all $x \in \Cantor$ and $b \in \Cantor$, there is some $y \in \Cantor$ which is $E$-equivalent to $x$ such that $x\oplus y\oplus a \geq_h b$.
\end{definition}

We will now show that for equivalence relations with all classes uncountable, Borel graphability implies the weak Borel pair-coding property. The proof largely follows the proof that $\Fwrel$ is not Borel graphable. In fact, that proof can be understood as consisting of two parts: first, a proof that Borel graphability implies the weak Borel pair-coding property and second, a proof that $\Fwrel$ does not have the weak Borel pair-coding property.

\begin{theorem}
Suppose $E$ is an analytic equivalence relation on $\Cantor$ whose equivalence classes are all uncountable. If $E$ is Borel graphable, then $E$ has the weak Borel pair-coding property.
\end{theorem}

\begin{proof}
Suppose that $G$ is a $\Delta^1_1(a)$ graphing of $E$. We claim that $E$ has the weak Borel pair-coding property as witnessed by the oracle $a$. Fix $x \in \Cantor$ and $b \in \Cantor$ and suppose for contradiction that there is no $y$ which is $E$-equivalent to $x$ such that $x\oplus y \oplus a \geq_h b$. We will show that the $E$-equivalence class of $x$ is countable, contradicting our assumption.

As mentioned above, the proof is similar to the proof that $\Fwrel$ is not Borel graphable. In particular, we will show that for every $y$, if there is a path from $x$ to $y$ in $G$ then $y \leq_h x\oplus a$, which is enough to finish the proof since there are only countably many reals hyperarithmetic in $x\oplus a$. We will prove this by induction on the length of the path.

The base case---when the length of the path is $0$, i.e.\ when $y = x$---is trivial. For the inductive case, let $z$ be the element immediately preceding $y$ on the shortest path from $x$ to $y$. By our inductive assumption, $z \leq_h x\oplus a$. Thus the set of neighbors of $z$ is $\Delta^1_1(x\oplus a)$.

By the relativized version of the effective perfect set theorem, either every neighbor of $z$ is hyperarithmetic in $x\oplus a$ (so, in particular, $y$ is hyperarithmetic in $x\oplus a$) or the set of neighbors of $z$ is uncountable. In the former case, we are done. In the latter case, the relativized version of Friedman's Conjecture implies that there is some neighbor $w$ of $z$ such that $w\oplus x\oplus a \geq_h b$, contradicting our assumption.
\end{proof}

We have now seen six different coding properties that an analytic equivalence relation can have: Borel witness coding, Borel coding, Borel pair-coding and their weak versions. The implications between these coding properties and with Borel graphability are summarized in the following diagram.

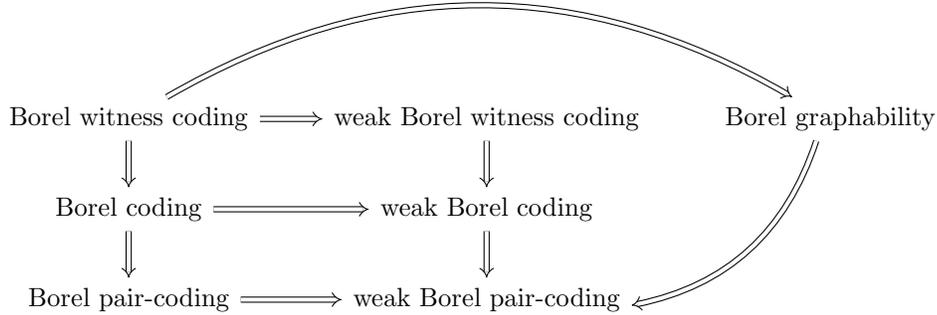
\begin{figure}[ht]
\centering
\begin{tikzcd}[arrows=Rightarrow]
\text{Borel witness coding} \arrow[r]\arrow[d]\arrow[rr, bend left] & \text{weak Borel witness coding} \arrow[d] & \text{Borel graphability} \arrow[ddl, bend left]\\
\text{Borel coding} \arrow[r]\arrow[d] & \text{weak Borel coding} \arrow[d] \\
\text{Borel pair-coding} \arrow[r] & \text{weak Borel pair-coding}
\end{tikzcd}
\caption{Implications between coding properties and graphability.}
\end{figure}

Note that the strongest of these coding properties, Borel witness coding, implies Borel graphability while the weakest, weak Borel pair-coding, is implied by Borel graphability (for equivalence relations whose classes are all uncountable). It would be nice to narrow the gap between these two results. One way to do this would be to strengthen the implications between coding and graphability, for example by showing that the Borel pair-coding property implies Borel graphability or that Borel graphability implies the Borel witness coding property. Another way would be to prove that some of the coding properties are equivalent to each other.

We will now put a limitation on these possibilities by separating some of the coding properties from each other, as well as from Borel graphability.

\begin{theorem}
\label{thm:separation}
There is an analytic equivalence relation $E$ such that $E$ is Borel graphable and has the Borel pair-coding property and the weak Borel witness coding property (and hence the weak versions of all of the coding properties) but does not have the Borel coding property.
\end{theorem}

The following diagram summarizes the consequences of this theorem. The equivalence relation $E$ from the theorem satisfies all circled properties and none of the properties which are not circled. Hence none of the circled properties imply any of the non-circled properties.

\begin{figure}[ht]
\centering
\begin{tikzcd}[arrows=Rightarrow,
    execute at end picture={
        \draw[draw, rounded corners, looseness = 20]
        (tikz@f@1-1-2.north east) -| (tikz@f@1-1-3.south east) .. controls (5, -.85) .. 
        (tikz@f@1-3-2.south east) -| (tikz@f@1-3-1.north west) -| (tikz@f@1-1-2.north west)
        -- (tikz@f@1-1-2.north east);
    }]
\text{Borel witness coding} \arrow[r]\arrow[d]\arrow[rr, bend left] & \text{weak Borel witness coding} \arrow[d] & \text{Borel graphability} \arrow[ddl, bend left]\\
\text{Borel coding} \arrow[r]\arrow[d] & \text{weak Borel coding} \arrow[d]  & {}\\
\text{Borel pair-coding} \arrow[r] & \text{weak Borel pair-coding} & {}
\end{tikzcd}
\caption{Non-implications between coding properties and graphability.}
\end{figure}
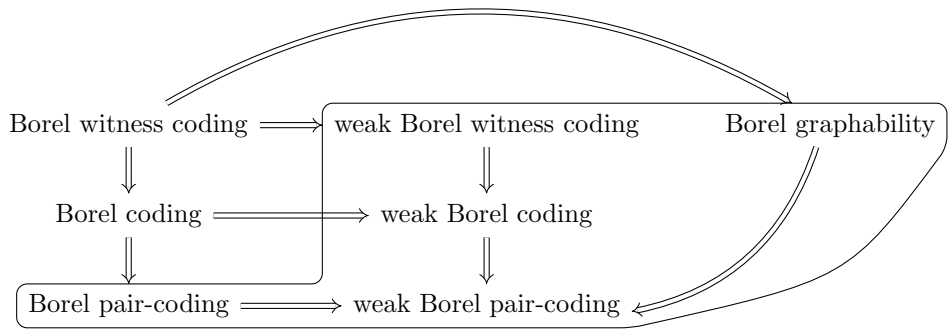

In order to define the equivalence relation $E$ that we will use to prove this theorem, we need to review some definitions from computability theory.

First, we need the notion of a jump hierarchy. Suppose $L$ is a linear order with domain $\N$ and $0$ is the least element of $L$. Given a real $a \in \Cantor$, a \term{jump hierarchy on $L$ which starts with $a$} is any $H \in 2^{\N\times\N}$ such that:
\begin{enumerate}
\item $H_{0} = a$
\item and for every $n \neq 0$, $H_n = \bigoplus_{m <_L n}H_m'$
\end{enumerate}
where $H_n$ denotes the $n^\text{th}$ column of $H$. Note that if $L$ is well-founded, then for every $a$, there is a unique jump hierarchy on $L$ starting from $a$. On the other hand, if $L$ is ill-founded then it is possible that there is no jump hierarchy on $L$ and it is possible that there are many different jump hierarchies on $L$. In the case where $L$ is well-founded, the unique jump hierarchy $H$ should be thought of as the $L^\text{th}$ iterate of the Turing jump of $a$, which justifies the following choice of notation.

\begin{notation}
If $L$ is well-founded then we will use $a^{(L)}$ to denote the unique jump hierarchy on $L$ which starts with $a$.
\end{notation}

One important point is that hyperarithmetic reducibility can be characterized in terms of jump hierarchies. In particular, $x \leq_h y$ if and only if there is some well-founded linear order $L$ which is computable from $y$ such that $x \leq_T y^{(L)}$. Similarly, if a function $f \colon \Cantor \to \Cantor$ is Borel then there is an oracle $a$ and a well-founded linear order $L$ which is computable from $a$ such that for all $x$, $f(x) \leq_T (a\oplus x)^{(L)}$.

Second, we need the following computability-theoretic notion of genericity. Given a set $S \subseteq 2^{< \N}$, a real $x \in \Cantor$ \term{meets or avoids $S$} if there is some finite initial segment $\sigma$ of $x$ such that either $\sigma$ is in $S$ or no string extending $\sigma$ is in $S$. Given reals $x$ and $a$, $x$ is \term{$a$-generic} if $x$ meets or avoids each subset of $2^{<\N}$ which is computably enumerable relative to $a$ (for a more in-depth introduction to genericity, see Section 2.24 of the book~\cite{downey2010algorithmic} of Downey and Hirschfeldt or the survey~\cite{jockusch1980degrees} by Jockusch). 

We will need the following fact about jump hierarchies of generic reals. It is well-known in computability theory; a version of it for finite iterates of the Turing jump can be found as Theorem 22.24.3 in~\cite{downey2010algorithmic}.

\begin{theorem}[Jockusch]
\label{thm:generic_jump}
For any real $a$, any well-founded linear order $L$ computable from $a$, and any $a^{(L)}$-generic $x$, $(a \oplus x)^{(L)} \equiv_T a^{(L)}\oplus x$.
\end{theorem}

It will also be useful to isolate the following corollary of this fact (which is also well-known in computability theory).

\begin{corollary}
For any reals $a$ and $x$, if $x$ is $a$-generic, then $a\oplus x$ does not compute $a'$.
\end{corollary}

\begin{proof}
Suppose for contradiction that $a \oplus x \geq_T a'$. It follows from the theorem above that $(a \oplus x)' \equiv_T a' \oplus x$. Hence we have
\[
    a \oplus x \geq_T a' \oplus x \equiv_T (a \oplus x)'.
\]
But this contradicts the fact that no real can compute its own jump.
\end{proof}

We are now ready to define the equivalence relation $E$ that we will use to prove Theorem~\ref{thm:separation}. Let $X$ be the Borel subset of $\lo\times\Cantor\times\Cantor$ consisting of those triples $(L, a, x)$ such that $L$ is a linear order with domain $\N$ which is computable from $a$ and which has $0$ as its least element. Define an equivalence relation $E$ on $X$ by setting $(L, a, x)$ and $(R, b, y)$ equivalent if $(L, a) = (R, b)$ and one of the following holds:
\begin{enumerate}
\item $L$ is ill-founded.
\item $L$ is well-founded and both $x$ and $y$ are $a^{(L)}$-generic.
\item $L$ is well-founded and neither $x$ nor $y$ is $a^{(L)}$-generic.
\end{enumerate}
In other words, corresponding to any pair $(L, a)$ there are either one or two equivalence classes. If $L$ is ill-founded then there is a single equivalence class and if $L$ is well-founded then there are two equivalence classes, one consisting of the reals which are $a^{(L)}$-generic and one consisting of the reals which are not $a^{(L)}$-generic.

\begin{proposition}
The equivalence relation $E$ is analytic.
\end{proposition}

\begin{proof}
Note that $(L, a, x)$ and $(L, a, y)$ are equivalent if and only if at least one of the following holds:
\begin{enumerate}
\item $L$ is ill-founded.
\item There is a jump hierarchy $H$ on $L$ which starts with $a$ such that both $x$ and $y$ are $H$-generic.
\item There is a jump hierarchy $H$ on $L$ which starts with $a$ such that neither $x$ nor $y$ are $H$-generic.
\end{enumerate}
Note that, unlike in the definition of $E$, in items (2) and (3) above we don't require that $L$ is well-founded. This is fine since if $L$ is ill-founded then $(L, a, x)$ and $(L, a, y)$ will be $E$-equivalent no matter what and so conditions (2) and (3) above don't matter. To finish, simply note that this characterization of $E$-equivalence is clearly $\Sigma^1_1$ and hence analytic.
\end{proof}

\begin{proposition}
\label{prop:pair_coding}
The equivalence relation $E$ has the Borel pair-coding property.
\end{proposition}

\begin{proof}
We will show that for any $(L, a, x)$ and any $b \in \Cantor$, there is some $y$ such that $(L, a, y)$ is $E$-equivalent to $(L, a, x)$ and $x\oplus y$ computes $b$. There are three cases to consider.

\medskip\noindent\textit{Case 1: $L$ is ill-founded.} This case is easy: we can just take $y = b$ since $(L, a, x)$ is $E$-equivalent to any triple of the form $(L, a, z)$.

\medskip\noindent\textit{Case 2: $L$ is well-founded and $x$ is not $a^{(L)}$-generic.} This case is also easy. Note that any real whose even bits are all $0$ is not $a^{(L)}$-generic. Thus we can take $y$ to be a real whose odd bits code $b$ and whose even bits are all $0$.

\medskip\noindent\textit{Case 3: $L$ is well-founded and $x$ is $a^{(L)}$-generic.} This is the only case that requires any real work. The idea is that we can build some $y$ which matches $x$ so closely that it inherits $x$'s genericity but which has information coded into the rare places where it disagrees with $x$.

We will begin by describing how to decode information from $x\oplus y$. Pick a sequence of numbers $n_0, n_1, n_2, \ldots$ as follows. First, let $n_0$ be the first position at which $x$ and $y$ disagree (i.e.\ at which $x(n_0) \neq y(n_0)$), then let $n_1$ be the first position strictly greater than $n_0 + 1$ at which $x$ and $y$ disagree, and so on. More generally, let $n_{k + 1}$ be the first position strictly greater than $n_k + 1$ at which $x$ and $y$ disagree. Now construct an element $z$ of $\Cantor$ by setting $z(k) = y(n_k + 1)$. We will think of this $z$ as the information encoded into $x\oplus y$. Recall that our goal is to encode $b$ into $x\oplus y$ and hence we need to construct $y$ in such a way that the real $z$ which we have just described is equal to $b$.

We will now explain how to actually build $y$. Let $S_0, S_1, S_2, \ldots$ enumerate the subsets of $2^{< \N}$ which are computably enumerable relative to $a^{(L)}$. We will build $y$ in a sequence of stages, at each stage determining some finite initial segment of $y$. On stage $n$, we will make that $y$ either meets or avoids $S_n$ and then code one more bit of $b$.

Suppose we have just finished stage $n$, let $\sigma$ be the initial segment of $y$ constructed so far and let $k$ be the length of $\sigma$; we will describe what to do on stage $n + 1$.

First, we will meet or avoid $S_{n + 1}$ by copying a long enough segment of $x$. For any $l$, let $x[k:k + l]$ denote the string given by taking the $l$ bits of $x$ starting at bit $k$ and ending at bit $k + l - 1$. We claim that there is some $l$ large enough that $\sigma\concat x[k:k + l]$ either meets or avoids $S_{n + 1}$. If not then $x$ fails to either meet or avoid the set
\[
  \{\tau_0\concat\tau_1 \mid |\tau_0| = k \text{ and } \sigma\concat\tau_1 \in S_{n + 1}\}
\]
which is clearly computable from $S_{n + 1}$ and hence c.e.\ in $a^{(L)}$. Pick some such large enough $l$ and set $\sigma' = \sigma\concat x[k:k + l]$.

Now we will code one more bit of $b$ into $y$. To do so, simply extend $\sigma'$ to $\sigma''$ by adding two more bits: first, a bit which disagrees with $x(k + l)$ and second, a bit which agrees with $b(n + 1)$. In other words, define $\sigma''$ to be
\[
\sigma'' = {\sigma'}\concat (1 - x(k + l))\concat b(n + 1)
\]
and finish stage $n + 1$ by setting $\sigma''$ to be an initial segment of $y$.
\end{proof}

\begin{proposition}
\label{prop:weak_coding}
The equivalence relation $E$ has the weak Borel witness coding property.
\end{proposition}

\begin{proof}
We will show that for any $(L, a, x)$ and any $b \in \Cantor$, there is some $y$ such that $(L, a, y)$ is $E$-equivalent to $(L, a, x)$ and a sufficient number of jumps of $a\oplus y$ computes $b$ as well as a witness to the $E$-equivalence of $(L, a, x)$ and $(L, a, y)$. As in the proof of Proposition~\ref{prop:pair_coding}, the only non-trivial part of the proof is the case when $L$ is well-founded and $x$ is $a^{(L)}$-generic. For convenience, we will assume that $b$ computes a witness to the fact that $x$ is $a^{(L)}$-generic.

Let $S_0, S_1, S_2, \ldots$ enumerate the subsets of $2^{< \N}$ which are c.e.\ in $a^{(L)}$. Build an $a^{(L)}$-generic real $y$ by alternately meeting or avoiding $S_n$'s and coding bits of $b$. More precisely, build $y$ in a sequence of stages, at each stage determining a finite initial segment of $y$.

Suppose we have just finished stage $n$ and let $\sigma$ denote the initial segment of $y$ formed so far. Let $\tau$ be the lexicographically least string extending $\sigma$ such that $\tau$ either meets or avoids $S_{n + 1}$. Finish stage $n + 1$ by setting $\tau\concat b(n + 1)$ to be an initial segment of $y$.

We now claim that $(a \oplus y)^{(L + 1)}$ computes both $b$ and a witness to the $E$-equivalence of $(L, a, x)$ and $(L, a, y)$. The main point is that $(a\oplus y)^{(L + 1)}$ computes $a^{(L + 1)}$ and hence can follow along the construction of $y$ in the previous paragraph, i.e.\ using $a^{(L + 1)}$ we can answer questions like ``what is the lexicographically least string that either meets or avoids $S_n$'' and this allows us to determine which bits of $y$ are coding bits of $b$. This also gives us a witness that $y$ is $a^{(L)}$-generic. Since $b$ was assumed to compute a witness that $x$ is $a^{(L)}$-generic, this gives us a witness to the $E$-equivalence of $(L, a, x)$ and $(L, a, y)$.
\end{proof}

\begin{proposition}
The equivalence relation $E$ is Borel graphable (with diameter $4$).
\end{proposition}

\begin{proof}
Define a graph $G$ by putting an edge between $(L, a, x)$ and $(L, a, y)$ if $x\oplus y$ computes a witness that $(L, a, x)$ and $(L, a, y)$ are equivalent. By a ``witness'' here we mean either of the following:
\begin{enumerate}
    \item A descending sequence in $L$.
    \item A jump hierarchy $H$ on $L$ which starts with $a$ such that either both $x$ and $y$ are $H$-generic or neither $x$ nor $y$ are $H$-generic.
\end{enumerate}
Note that, given $H$ and $x$, it is Borel to check whether $x$ is $H$-generic.

It is clear that $G$ is a subset of $E$, so we just need to show that given any $(L, a, x)$, the equivalence class of $(L, a, x)$ is connected in $G$. As in the proofs of Propositions~\ref{prop:pair_coding} and~\ref{prop:weak_coding} above, this is easy to do when $L$ is ill-founded or when $L$ is well-founded and $x$ is not $a^{(L)}$-generic. Thus we just need to deal with the case that $L$ is well-founded and $x$ is $a^{(L)}$-generic.

We will need the following definition and fact about genericity. Given any real $a$, reals $x$ and $y$ are \term{mutually $a$-generic} if $x\oplus y$ is $a$-generic. It is well-known (and not hard to prove) that if $x$ and $y$ are both $a$-generic, then there is some $z$ such that $z$ is $a$-generic and mutually $a$-generic with respect to both $x$ and $y$.

Now let $L$ and $a$ be as above. We will prove that if $x$ and $y$ are mutually $a^{(L)}$-generic then $(L, a, x)$ and $(L, a, y)$ are connected in $G$ by a path of length $2$. This is enough to finish the proof since given an arbitrary $x$ and $y$ which are $a^{(L)}$-generic, we can find some $z$ which is $a^{(L)}$-generic and mutually $a^{(L)}$-generic with respect to both $x$ and $y$. Thus $(L, a, z)$ is connected in $G$ to both $(L, a, x)$ and $(L, a, y)$ by paths of length $2$, implying that $(L, a, x)$ and $(L, a, y)$ are connected by a path of length $4$.

We now assume that $x$ and $y$ are mutually $a^{(L)}$-generic. We will construct an $a^{(L)}$-generic real $w$ such that $x\oplus w$ and $y\oplus w$ can both compute $a^{(L)}$. This gives us that $(L, a, w)$ is adjacent in $G$ to both $(L, a, x)$ and $(L, a, y)$, thus finishing the proof.

The construction of $w$ is very similar to the proof of Proposition~\ref{prop:pair_coding}. That is, we will basically take $w = x\oplus y$, but occasionally change some bits to encode information; to encode information so that $x$ can decode it, we will alter even bits of $w$ and to encode information so that $y$ can decode it, we will alter odd bits. Genericity of $w$ comes from the fact that $x$ and $y$ are mutually generic.

We now describe the construction in more detail. Suppose that by the end of stage $n$, we have constructed an initial segment $\sigma$ of $w$. We first extend $\sigma$ to $\sigma'$ by matching the bits of $x\oplus y$ in order to meet or avoid the next $a^{(L)}$-c.e.\ set of strings. This is always possible because $x\oplus y$ is $a^{(L)}$-generic. Next, we add bits onto $\sigma'$ which disagree with $x$ and $y$ and then add bits to encode the next bit of the information we are trying to encode (namely $a^{(L)}$). More precisely, suppose that the next bit of information to encode is $b$ and that $\sigma'$ has length $2k$. Then we extend $\sigma'$ to
\[
    (\sigma')\concat (1 - x(k)) \concat (1 - y(k)) \concat b\concat b
\]
and declare this to be an initial segment of $w$, thus finishing stage $n + 1$.

Given $w$ constructed in this way, $x\oplus w$ can compute the encoded information in much the same manner as in the proof of Proposition~\ref{prop:pair_coding}, except that when using $x$ to decode information, we only look at the even bits of $w$ and when using $y$, we only look at the odd bits.\footnote{Note that in our construction of $w$ in the paragraph above, we appended the bit $b$ twice so that it is encoded into both the even and the odd bits of $w$.}
\end{proof}

\begin{proposition}
The equivalence relation $E$ does not have the Borel coding property.
\end{proposition}

\begin{proof}
Suppose for contradiction that $f \colon X \to \Cantor$ is a Borel function witnessing that $E$ has the Borel coding property. Since $f$ is Borel, there is some real $a$ and some well-order $L$ computable from $a$ such that for all $(R, b, x)$ in $X$, $f(R, b, x)$ is computable from $(a\oplus R\oplus b\oplus x)^{(L)}$. Let $x$ be $a^{(L)}$-generic. We claim that for any $(L, a, y)$ which is $E$-equivalent to $(L, a, x)$, $f(L, a, x)$ does not compute $a^{(L + 1)}$, contradicting the assumption that $f$ witnesses the Borel coding property for $E$.

To prove our claim, note that if $(L, a, y)$ is $E$-equivalent to $(L, a, x)$ then $y$ is $a^{(L)}$-generic. Thus we can compute
\begin{align*}
f(L, a, y) &\leq_T (a\oplus L \oplus a\oplus y)^{(L)} &\text{(by our choice of $a$ and $L$)}\\
           &\leq_T (a\oplus y)^{(L)} &\text{(since $a$ computes $L$)}\\
           &\leq_T a^{(L)}\oplus y &\text{(by Theorem~\ref{thm:generic_jump}).}
\end{align*}
And since $y$ is $a^{(L)}$-generic, $a^{(L)}\oplus y$ does not compute $a^{(L + 1)} = (a^{(L)})'$.
\end{proof}

All this leaves open the question of whether Borel graphability of analytic equivalence relations with all classes uncountable is equivalent to some coding property. We have seen above that it cannot be equivalent to either the Borel witness coding property or the Borel coding property, but perhaps it could be equivalent to the Borel pair-coding property or to one of the weak coding properties.


\section{Hypergraphability}
\label{sec:hypergraph}

In this section we will consider a variation of Borel graphability defined in terms of hypergraphs rather than graphs. Actually, there are at least two reasonable variations of this sort. The first, which we will call \term{Borel hypergraphability}, is the most straightforward generalization of Borel graphability to hypergraphs. The second, which we will call \term{strong Borel hypergraphability}, can be seen as a generalization of the notion of being Borel graphable by a graph of finite diameter.

\begin{definition}
\sloppy An analytic equivalence relation $E$ on $X$ is \term{Borel {$n$-hypergraphable}} if there is a Borel hypergraph $G$ on $X$ such that:
\begin{enumerate}
\item Each hyperedge of $G$ has at most $n$ points.
\item If $\{x_1, \ldots, x_n\}$ is a hyperedge of $G$ then $x_1, \ldots, x_n$ are all in the same $E$-equivalence class.
\item If $x E y$ then there is some finite sequence of hyperedges $X_1, \ldots, X_m$ such that $x$ is in $X_1$, $y$ is in $X_m$ and each $X_i$ has nonempty intersection with $X_{i + 1}$.
\end{enumerate}
Also $E$ is \term{strongly Borel $n$-hypergraphable} if there is some Borel hypergraph $G$ which satisfies the above conditions, but with the third condition strengthened to the following: if $x E y$ then there is some hyperedge $X$ such that $x$ and $y$ are both in $X$.
\end{definition}

\begin{definition}
An analytic equivalence relation $E$ on $X$ is \term{Borel hypergraphable} if it is Borel $n$-hypergraphable for some $n$ and is \term{strongly Borel hypergraphable} if it is strongly Borel $n$-hypergraphable for some $n$.
\end{definition}

Note that $E$ is Borel $2$-hypergraphable if and only if it is Borel graphable and strongly Borel $2$-hypergraphable if and only if it is Borel. Furthermore, if $E$ is Borel graphable with diameter $n$, then it is strongly Borel $(n + 1)$-hypergraphable, which makes precise our comment above that strong Borel hypergraphability is a generalization of the property of being Borel graphable with finite diameter. 

Also note that Borel hypergraphability shares some basic properties with Borel graphability. In particular, if $E$ is Borel hypergraphable, then it is analytic and if all equivalence classes of $E$ are countable, then $E$ is Borel hypergraphable if and only if $E$ is Borel. This second point shows that there are analytic equivalence relations which are not Borel hypergraphable; see Proposition~\ref{prop:firstnonexample} for the details.

Our main reason for introducing Borel hypergraphability is that there are several questions about Borel graphability which we were unable to resolve, but for which we can answer the analogous question about Borel hypergraphability. In particular, we can show that every Polish group action induces a Borel hypergraphable equivalence relation and that the Borel pair-coding property implies Borel hypergraphability. 

\begin{theorem}
\label{thm:polish}
If $E$ is the orbit equivalence relation of a Polish group action $\Gamma \acts X$, then $E$ is strongly Borel $3$-hypergraphable.
\end{theorem}

\begin{proof}
As in the proof of Theorem \ref{thm:coding_to_graphable}, we can assume that $\Gamma=X=\Cantor$ and the group operations and the action are Borel. We can also assume that all orbits are uncountable, so by Corollary \ref{cor:polish_group_pair_coding}, we have Borel functions $e \colon X\times \Cantor \to X$, $d\colon X\times X \to \Cantor$, and $g\colon X\times X \to \Gamma$ such that for any $x \in X$ and $a \in \Cantor$, if $y = e(x, a)$ then the following hold:
\begin{enumerate}
\item $d(x, y) = a$ and
\item  $y$ is in the same orbit as $x$ as witnessed by $g(x, y)$---i.e.\ $g(x, y)\cdot x = y$.
\end{enumerate}

Now define a $3$-hypergraph on $X$ as follows. Given $\{x, y, z\}$, let $A_{x, y, z}$ denote the smallest subset of $X\cup \Gamma\cup \Cantor$ which contains $x, y$ and $z$ and which is closed under the group operations of $\Gamma$, the group action of $\Gamma$ on $X$, the functions $e, d, g$, finite joins (in the sense of computability theory) and Turing reductions. Declare $\{x, y, z\}$ to be a hyperedge if this set $A_{x, y, z}$ contains group elements taking each of $x, y,$ and $z$ to the others.

Now suppose we are given $x$ and $y$ in the same orbit; we will show there is some hyperedge that contains both of them. Let $a \in \Gamma$ be such that $a\cdot x = y$. Set $z = e(x, a)$. We claim that $\{x, y, z\}$ forms a hyperedge. It is enough to show that we can obtain group elements taking $x$ to $z$ and $x$ to $y$ (the rest follows by multiplying those elements together). For the former, simply note that $g(x, z)$ is a group element taking $x$ to $z$. For the latter, note that $d(x, z)$ is a group element taking $x$ to $y$.
\end{proof}

\begin{theorem}
If $E$ is an analytic equivalence relation with the Borel pair-coding property, then $E$ is Borel $3$-hypergraphable and strongly Borel $4$-hypergraphable.
\end{theorem}

\begin{proof}
We will first prove that $E$ is Borel $3$-hypergraphable and then observe how to modify the proof to get strong $4$-hypergraphability.

For convenience we assume the domain of $E$ is $\Cantor$. Let $f \colon \Cantor\times\Cantor \to \Cantor$ witness the Borel pair-coding property for $E$. The hypergraph is defined by declaring that $\{x_1, x_2, x_3\}$ forms a hyperedge if together $x_1, x_2, x_3, f(x_1, x_2), f(x_1, x_3),$ and $f(x_2, x_3)$ compute sufficient evidence that all of $x_1, x_2,$ and $x_3$ are $E$-equivalent. Slightly more formally, here's what we mean by ``sufficient evidence.'' Consider the set of all reals computable from $x_1, \ldots, f(x_1, x_3)$. Any time two of the reals in this set are $E$-equivalent, as witnessed by another of the reals in the set, put an edge between them. Then declare $\{x_1, x_2, x_3\}$ to be a hyperedge if they are in the same connected component of the resulting graph.

Now suppose that $x$ and $y$ are $E$-equivalent. Since $E$ has the Borel pair-coding property, every equivalence class of $E$ must be uncountable. In particular, the equivalence class of $x$ and $y$ must be uncountable. An examination of the proof of the effective perfect set theorem then tells us that there is a perfect, binary-branching tree $T$ on $2\times\N$ such that the following holds:
\begin{enumerate}
    \item For any $(z, a) \in [T]$, $x E z$ and $a$ is a witness to this fact.
    \item For any distinct $(z_1, a_1), (z_2, a_2) \in [T]$, $z_1 \neq z_2$.
\end{enumerate}
The reason that the second item is useful is that it ensures that if $z$ is the first coordinate of some element $(z, a)$ of $[T]$ then $z\oplus T$ computes $a$---in other words, $z \oplus T$ computes a witness to the $E$-equivalence of $z$ and $x$. Also note that since $T$ is a perfect, binary branching tree, for any real $b$, we can find some path $(z, a) \in [T]$ such that $(z, a)\oplus T$ computes $b$ (because we can code the bits of $b$ into decisions about whether to turn left or right at each branch point of $T$).

Since $f$ witnesses that $E$ has the Borel pair-coding property, we can find $z_1$ and $z_2$ such that $x E z_1$, $y E z_2$ and $f(x, z_1), f(y, z_2)$ both compute $T$. Now let $(w, a)$ be some element of $[T]$ such that $(w, a) \oplus T$ computes $x$ along with witnesses to the $E$-equivalence of $x, y, z_1$ and $z_2$. 

We now claim that $\{x, z_1, w\}$ and $\{y, z_2, w\}$ are both hyperedges. Here's why $\{x, z_1, w\}$ is a hyperedge:
\begin{enumerate}
    \item From $f(x, z_1)$, we can compute $T$.
    \item From $w$ and $T$, we can compute $a$, which is a witness to the $E$-equivalence of $x$ and $w$. 
    \item From $(w, a)$ and $T$, we can compute witnesses to the $E$-equivalence of $x$ and $z_1$.
\end{enumerate}
The proof for $\{y, z_2, w\}$ is similar except that instead of being able to directly compute a witness that $y$ and $w$ are $E$-equivalent, we have to compute $x$ and witnesses that $w$ and $x$ as well as $x$ and $y$ are $E$-equivalent.

To prove strong $4$-hypergraphability, define the hypergraph in the same way, except allowing four reals in each hyperedge instead of three. Then note that if $x, y, z_1,$ and $w$ are as above then $\{x, y, z_1, w\}$ forms a hyperedge in this new hypergraph.
\end{proof}

Based on the results of this section up to this point, one might think that Borel graphability and Borel hypergraphability are essentially equivalent and the fact that we are able to prove more for Borel hypergraphability only reflects our lack of techniques for Borel graphability. However, by considering the relations $\Fw$ and $\Fwrel$ from Section~\ref{sec:admiss}, we can see that this is not the case. In particular, we saw earlier that the Borel graphability of $\Fw$ is independent of $\ZFC$ and that $\Fwrel$ is never Borel graphable, but we will now show that both of them are Borel hypergraphable (provably in $\ZFC$).

\begin{theorem}
Both $\Fw$ and $\Fwrel$ are Borel $3$-hypergraphable and strongly Borel $4$-hypergraphable.
\end{theorem}

\begin{proof}
We will only prove this for $\Fw$. The result for $\Fwrel$ follows by relativizing the proof. We will first prove that $\Fw$ is Borel $3$-hypergraphable and then observe how to modify the proof to get strong $4$-hypergraphability.

The hypergraph is defined by declaring that $\{x_1, x_2, x_3\}$ forms a hyperedge if together $x_1, x_2$ and $x_3$ compute witnesses that each one is $\Fw$-equivalent to the others. 

Now suppose we are given $x$ and $y$ such that $\omega_1^x = \omega_1^y$. The main idea is to follow the proof of Lemma~\ref{lemma:omega_1_L} (which was the main step in proving that $\Fw$ is Borel graphable if there is a non-constructible real) but to replace the reliance on the parameter in the definition of the graph by adding the parameter to the hyperedge.

Let $M$ be a countable transitive model of (a sufficiently large finite fragment of) $\ZFC$ such that $M$ contains $x$ and let $a$ be a real which is not in $M$ and which satisfies $\omega_1^a = \omega_1^x$ (such an $a$ must exist since there are uncountably many reals $z$ with $\omega_1^z = \omega_1^x$ and $M$ is countable). We now build a real $g$ which is Kumabe-Slaman generic over $M$ with the aim of ensuring that $\{x, a, x\oplus g\}$ and $\{y, a, x\oplus g\}$ are both hyperedges.

Since $g$ is Kumabe-Slaman generic over $M$ and $M$ contains $x$, we have that $\omega_1^{x\oplus g} = \omega_1^x$. And by following the proof of Lemma~\ref{lemma:omega_1_L}, we can build $g$ such that $g\oplus a$ computes witnesses to the $\Fw$-equivalence of $x, y, a,$ and $x\oplus g$. Thus $\{x, a, x\oplus g\}$ and $\{y, a, x\oplus g\}$ are hyperedges, as promised.

To prove strong $4$-hypergraphability, define a new hypergraph by declaring $\{x_1, x_2, x_3, x_4\}$ to form a hyperedge if together they compute witnesses that each one is $\Fw$-equivalent to the others. Then note that if $x, y, a,$ and $g$ are as above then $\{x, y, a, x\oplus g\}$ forms a hyperedge in this new hypergraph.
\end{proof}

\begin{corollary}
There is an equivalence relation which is Borel hypergraphable but not Borel graphable.
\end{corollary}

The results of this section give the impression that most analytic equivalence relations are Borel hypergraphable. We know that this is not true for all analytic equivalence relations because, as we mentioned earlier, an equivalence relation whose equivalence classes are all countable is Borel hypergraphable if and only if it is Borel. However, to the best of our knowledge it could be true for analytic equivalence relations with no countable equivalence classes.

\begin{question}
Is there an analytic equivalence relations whose equivalence classes are all uncountable but which is not Borel hypergraphable?
\end{question}


\section{Computably graphable equivalence relations}
\label{sec:ceers}

There is an analogy between analytic equivalence relations on $\Cantor$ and computably enumerable equivalence relations on $\N$ (which are usually referred to as ``ceers''). Extending this analogy to Borel graphability, we have the following definition.

\begin{definition}
A ceer $E$ is \term{computably graphable} if $E$ is the connectedness relation of a computable graph $G$ on $\N$. $G$ is called a \term{computable graphing} of $E$.
\end{definition}

It is not hard to find examples of ceers which are not computably graphable. For example, suppose $A$ is a computably enumerable, non-computable set and let $E$ be the equivalence relation given by:
\[
    a E b \iff a = b \text{ or } \exists n\, (a = 2n,\, b = 2n + 1 \text{ (or vice-versa) and } n \in A).
\]
Note that all equivalence classes of $E$ have size one or two. From this, it is easy to conclude that $E$ is computably graphable if and only if it is computable. But if $E$ is computable then so is $A$, hence $E$ is not computably graphable.

It is easy to modify the example above to give a non-computably graphable ceer with equivalence classes of any finite size. However, if we restrict to equivalence relations where all equivalence classes are infinite, then every ceer is computably graphable.\footnote{A related result was independently obtained by Uri Andrews and coauthors (unpublished).} This is in contrast to the situation with analytic equivalence relations on $\Cantor$, where there are non-Borel graphable analytic equivalence relations, all of whose equivalence classes are uncountable.

\begin{theorem}
Suppose $E$ is a ceer whose equivalence classes are all infinite. Then $E$ is computably graphable with diameter $2$.
\end{theorem}

\begin{proof}
Since $E$ is computably enumerable, there is a computable procedure, which, for each $x \in \N$, enumerates the elements of $x$'s equivalence class (without duplicates). More precisely, there is a computable function $f \colon \N \times \N \to \N$ such that $f(x, 0), f(x, 1), \ldots$ enumerates the equivalence class of $x$ without duplicates. Furthermore, since each equivalence class is infinite, $f$ is a total function. Note, by the way, that even if $x$ and $y$ are equivalent, $f(x, -)$ and $f(y, -)$ may enumerate the equivalence class in different orders.

Now define a computable graph as follows. Suppose we are given $x < y$ and want to decide if they are adjacent. Then for every $z < y$, we enumerate the first $y$ elements of $z$'s equivalence class---i.e.\ we look at $f(z, 0), f(z, 1), \ldots, f(z, y)$. This gives us a finite subset of $E$ which we can think of as a finite graph (where the vertices are all $z$'s below $y$ along with anything else that happened to be enumerated and there is an edge between $u$ and $v$ if $v$ showed up in $u$'s enumeration or vice-versa). We say that $x$ and $y$ are adjacent if there is path between them in this graph. Note that if $x$ and $y$ are adjacent then they must be $E$-equivalent because this finite graph is a subset of $E$.

Now fix $x, y$ which are $E$-equivalent. We will show that there is some $z$ adjacent to both. Let $N$ be large enough that $y$ has shown up by stage $N$ of $x$'s enumeration (i.e. for some $n < N$, $f(x, n) = y$) and $x, y < N$. Let $z$ be the first element of $x$'s equivalence class enumerated by $x$ which is larger than $N$ (i.e. $z = f(x, n)$ for the least $n$ such that $f(x, n) > N$). 

We claim that $z$ is adjacent to both $x$ and $y$. First note that $z$ must show up by stage $z$ in $x$'s enumeration. This is because everything that $x$ enumerates before $z$ is $\leq N$ (by definition of $z$) and hence there are at most $N + 1$ things enumerated before $z$. This shows that $x$ and $z$ are adjacent.

Now let's prove that $z$ and $y$ are adjacent. Since $N < z$, if we start enumerating elements of $x$'s equivalence class then $y$ will show up before stage $z$ (as will $z$ itself, as noted above). Thus if we enumerate $z$ elements of each equivalence class of numbers below $z$ then we will see that $y$ and $z$ are equivalent.
\end{proof}

\appendix

\section{Kumabe-Slaman forcing}

In this appendix, we will provide an exposition of Kumabe-Slaman forcing and prove the facts about it that we use in Section~\ref{sec:admiss}. 

Kumabe-Slaman forcing was invented by Kumabe and Slaman in the early 1990s in order to prove generalizations of the Posner-Robinson theorem in computability theory. However, it only appeared in print several years later, in a paper by Shore and Slaman~\cite{shore1999defining}. Prior expositions of the forcing notion can be found in the paper by Shore and Slaman as well as papers by Jananthan and Simpson~\cite{jananthan2021turing} and by Reimann and Slaman~\cite{reimann2022effective}. 

We have two main reasons for presenting another exposition in this appendix. First, our presentation is somewhat different than that of the papers mentioned above (though it was influenced by informal explanations given by Slaman) and we feel that our presentation may be easier for those outside of computability theory to understand. Second, in Section~\ref{sec:admiss} we used one result about Kumabe-Slaman forcing which has not been published before. Namely, the fact that Kumabe-Slaman forcing preserves admissibility of countable ordinals. This fact was discovered in joint work by the third author of this paper and Benjamin Siskind. Since the proof has not yet been made publicly available, we have decided to include it in this appendix and to include an account of Kumabe-Slaman forcing in order to keep the appendix self-contained.

\subsection{The end result of Kumabe-Slaman forcing}\label{forcing}

The easiest way to understand Kumabe-Slaman forcing is to first understand what sort of object it is meant to produce. So before we give a precise definition of the forcing notion, we will explain that.

Essentially, a Kumabe-Slaman generic is a partial labelling of the tree $2^{<\N}$ with label set $\{0, 1\}$. Informally, this means that every node in the tree $2^{<\N}$ is either given a label of $0$ or $1$ or given no label at all. Formally, we can define such a labelling as a function $g \colon 2^{<\N} \to \{0, 1, \bot\}$ where $g(\sigma) = \bot$ indicates that the string $\sigma$ has not been given a label. 

We can use any such labelling $g$ to define a function $\Cantor \to 2^{\leq \N}$ as follows. Given any $x \in \Cantor$, we can think of $x$ as describing a path through $2^{<\N}$. Imagine that we follow this path and each time we encounter a node to which $g$ has assigned a label, we write that label down. In the end, we will be left with a sequence of $0$s and $1$s which could be either infinite or finite (depending on whether or not there are infinitely many labeled nodes along the path described by $x$). In other words, we will be left with an element of $2^{\leq \N}$. We can then define a function by mapping each $x$ to the element of $2^{\leq \N}$ resulting from this process. This is pictured below.

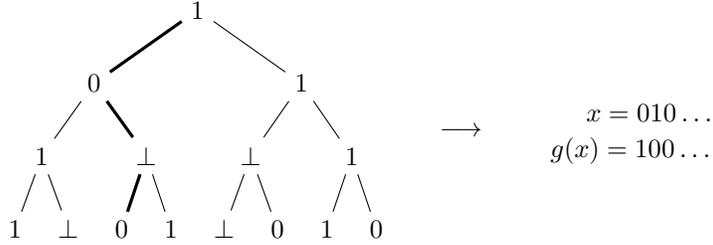
\begin{figure}[ht]
\centering
\begin{forest}
[1
  [0, edge={very thick}
    [1
        [1]
        [$\bot$]
    ]
    [$\bot$, edge={very thick}
        [0, edge={very thick}]
        [1]
    ] 
  ]
  [1
    [$\bot$
        [$\bot$]
        [0]
    ]
    [1
        [1]
        [0]
    ]
  ]
]
\node at (3.5, -1.5) {$\longrightarrow$};
\node at (6, -1.25) {$x = 010\ldots$};
\node at (5.75, -1.75) {$g(x) = 100\ldots$};
\end{forest}
\caption{A depiction of computing $g(x)$. The labelled tree on the left represents $g$. The bold edges in the tree denote the path through the tree determined by the bits of $x$. The sequence $g(x)$ is given by writing down the labels along this path.}
\end{figure}

To keep our notation simple, we will often conflate a labelling $g$ with the associated function that we have just described. In particular, we will use $g(x)$ to denote the output of this function on input $x$. Note that the function $x \mapsto g(x)$ is clearly computable relative to $g$.

We can now explain the key property of Kumabe-Slaman forcing. Suppose that $M$ is a countable transitive model of $\ZFC$ and $x$ is a real not in $M$. Then using Kumabe-Slaman forcing over $M$, one can construct a labelling $g$ such that $g(x)$ encodes essentially any information---including information about $g$ itself. Note that this includes information about $g$ which is highly non-computable relative to $g$, e.g.\ information about the hyperjump of $g$.

\subsection{The forcing notion}

We will now actually define the Kumabe-Slaman forcing notion.

\begin{definition}
A \term{condition for Kumabe-Slaman forcing} consists of a pair $p = (g_p, X_p)$ where $g_p$ is a partial labelling of a finite initial segment of $2^{<\N}$---formally a function $g_p\colon 2^{\leq n} \to \{0, 1, \bot\}$ for some $n$---and $X_p$ is a finite subset of $\Cantor$.
\end{definition}

Given a condition $p = (g_p, X_p)$, the set $X_p$ should be thought of as a set of ``forbidden paths''---paths along which we are not allowed to add any more labels later in the forcing.

We now note one somewhat non-obvious feature of this definition. Suppose that $p = (g_p, X_p)$ is a condition for Kumabe-Slaman forcing with $g_p \colon 2^{\leq n} \to \{0, 1, \bot\}$. Then there are two types of nodes to which $g_p$ does not assign a label: first, nodes $\sigma$ for which $|\sigma| \leq n$ and $g_p(\sigma) = \bot$ and second, nodes $\sigma$ for which $|\sigma| > n$ (so that $g_p(\sigma)$ is not even defined). One can think of the first type of node as a node which $p$ has \emph{explicitly decided} to not assign a label to and the second type as a node to which $p$ has not yet decided what label (if any) to assign. The point is that a condition extending $p$ can only add a label to the second type of node, not the first.

\begin{definition}
Given two conditions $p = (g_p, X_p)$ and $q = (g_q, X_q)$, $q$ \term{extends} $p$, written $q \leq p$, if all of the following hold.
\begin{enumerate}
\item $g_p \subseteq g_q$. This means that if $p$ assigns a label to some node then $q$ assigns the same label to that node. It also means that if $p$ has explicitly decided not to assign a label to a node (i.e. if $g_p(\sigma) = \bot$) then $q$ will also explicitly not assign a label to that node (i.e. $g_q(\sigma) = \bot$).
\item $X_p \subseteq X_q$. In other words, every path forbidden by $p$ is still forbidden by $q$.
\item For any $x \in X_p$ and any $\sigma \leq x$, if $g_p(\sigma)$ is not defined then $g_q$ does not assign a label to $\sigma$ (i.e.\ either $g_q(\sigma) = \bot$ or $g_q(\sigma)$ is not defined). In other words, $q$ does not add any new labels along paths in $X_p$.
\end{enumerate}
\end{definition}

Technically, a generic filter for Kumabe-Slaman forcing is a collection of Kumabe-Slaman conditions. However, we said at the beginning that Kumabe-Slaman forcing produces a partial labelling of $2^{< \N}$. Here's what this means. Given a generic filter $G$, we can define a partial labelling $g$ by setting
\[
  g = \bigcup_{p \in G}g_p.
\]
In other words, $g$ assigns a label to a node $\sigma \in 2^{< \N}$ if some condition in $G$ assigns that label to $\sigma$. Note that when defining $g$, we have remembered only the labels assigned by conditions in $G$ and have ``forgotten'' the forbidden paths.

We will typically engage in a slight abuse of terminology by referring to $g$ as a Kumabe-Slaman generic (even though it is really $G$ that is the generic). However, when we need to be precise, we will refer to $g$ as the \term{partial labelling associated to $G$}.

\subsection{The key lemma}

We will now prove the key lemma about Kumabe-Slaman forcing. Our statement (and proof) of the lemma is taken from Jananthan and Simpson~\cite[Lemma 3.7]{jananthan2021turing}, who were in turn inspired by a lemma of Reimann and Slaman~\cite[Lemma 3.10]{reimann2022effective}. The lemma, while somewhat technical looking, will be used to derive several important properties of Kumabe-Slaman forcing, including the ``key property of Kumabe-Slaman forcing'' that we discussed above.

For one of our uses of this lemma, it will be important to state it with a fairly high level of generality. In particular, we will state it in terms of Kumabe-Slaman forcing over (possibly ill-founded) $\omega$-models of $\ZF$, rather than just over transitive models (which are always well-founded). We will also need to pay some attention to what sorts of things are absolute between such models.

Recall that an $\omega$-model of $\ZF$ (or of other similar theories) is a model $M$ such that the natural numbers in $M$ are isomorphic to the true natural numbers. In practice we will typically ignore the words ``isomorphic to'' and just assume that the natural numbers in $M$ literally are the true natural numbers.

Note that for such a model $M$, many concepts are the same when interpreted inside $M$ as they are ``in the real world.'' For example, suppose $M \models x \in \Cantor$---i.e.\ $x$ is an element of $M$ which $M$ believes to be an element of $\Cantor$. Then $x$ really does correspond to an element of $\Cantor$. More precisely, we can define a real $\widetilde{x}$ by
\[
  \widetilde{x}(n) = 0 \iff M \models x(n) = 0
\]
and any arithmetic property of $\tilde{x}$ will hold if and only if $M$ believes the property holds of $x$. From this point on, we will not distinguish between the element $x$ of $M$ and the real $\widetilde{x}$ which corresponds to it.

Note that this means that if $p$ is an element of $M$ which $M$ believes to be a condition for Kumabe-Slaman forcing then $p$ really is a condition for Kumabe-Slaman forcing. This is because such a condition just consists of a finite object, along with a finite number of reals. $M$ is also correct about the ordering on conditions. We will often use the phrase ``Kumabe-Slaman forcing over $M$'' to mean the Kumabe-Slaman forcing poset as interpreted in $M$; the above reasoning shows that this is really just a sub-poset of the usual Kumabe-Slaman forcing poset.

We can now state the key lemma.

\begin{lemma}[Key Lemma of Kumabe-Slaman forcing]
\label{lemma:key}
Suppose we have:
\begin{itemize}
\item $M$, an $\omega$-model of $\ZF$.
\item $D \in M$, a dense set for Kumabe-Slaman forcing over $M$.
\item $p \in M$, a condition for Kumabe-Slaman forcing over $M$.
\item $X \subseteq \Cantor$, a finite set of reals which are not necessarily in $M$.
\end{itemize}
Then there is a condition $q \leq p$ in $M$ such that $q$ is in $D$ and $q$ adds no new labels along paths in $X$.
\end{lemma}

\begin{proof}
Suppose for contradiction that there is no such condition $q$---in other words, every condition $q \in M$ which extends $p$ and is in $D$ assigns at least one new label along some path in $X$. Here's the main idea of the proof. We first show that there is a finite set $Y$ of reals in $M$ which imitates $X$ in the sense that every condition $q \in M$ which extends $p$ and is in $D$ assigns at least one new label along some path in $Y$. We can then extend the condition $p$ to a condition $p' \in M$ by adding these imitators as forbidden paths. Since the $Y$ imitates $X$ in the sense just explained, it is impossible for any $q$ extending $p'$ to be in $D$. But this contradicts the density of $D$.

Let $x_1,\ldots, x_n$ enumerate the elements of $X$. Say that an $n$-tuple of strings $(\sigma_1,\ldots, \sigma_n) \in (2^{<\N})^n$ is \term{essential} if for every condition $q \leq p$ in $M$, if $q \in D$ then $q$ assigns a new label to some node which is compatible with one of the $\sigma_i$'s. In other words, $q$ assigns a new label to a node $\tau$ such that for some $i$, either $\sigma_i \leq \tau$ or $\tau \leq \sigma_i$.

Now define a tree $T \subseteq (2^{<\N})^n$ as follows. The nodes in $T$ consist of tuples of strings $(\sigma_1,\ldots,\sigma_n)$ such that all the $\sigma_i$'s are of the same length and the tuple is essential. It is not hard to check that $T$ is closed under taking component-wise initial segments and thus it is really a tree. Furthermore, it is clear from the definition of $T$ that the tuple $(x_1,\ldots, x_n)$ constitutes an infinite path through $T$ and thus $T$ is infinite.

Note that since $M$ is an $\omega$-model, the infinitude of $T$ is absolute between $M$ and $V$. Thus since $T$ is infinite (as we have just shown), $M$ believes it is infinite. And since $T$ is finitely branching (at most $2^n$-branching to be precise), $M$ believes that there is an infinite path through $T$.

For the remainder of the proof, we will officially reason inside $M$. However, since $M$ is an $\omega$-model, everything we will talk about is actually absolute between $M$ and $V$ and so this is just a formality.

Let $(y_1,\ldots,y_n)$ be an infinite path through $T$ and set $Y = \{y_1,\ldots,y_n\}$. This set $Y$ is the set of ``imitators'' that we mentioned at the beginning of the proof (recall that we are reasoning in $M$ and so we have $Y \in M$). We will now prove this formally.

To be precise, we must show that if a condition $q \in M$ extends $p$ and is in $D$, then $q$ adds at least one new label along some path in $Y$. To see why, suppose we have such a $q$. Since $g_q$ is finite, there is some number $N$ such that every node to which $q$ assigns a label has height less than $N$. Since $(y_1,\ldots, y_n)$ is a path through $T$, $(y_1\restriction N, \ldots, y_n\restriction N)$ is essential. In particular, $q$ must assign a new label to some node which is compatible with one of $y_1\restriction N, \ldots, y_n\restriction N$. But since $q$ assigns labels only to nodes of height less than $N$, this means that $q$ actually must assign a new label to a node which is an initial segment of one of $y_1\restriction N, \ldots, y_n\restriction N$ (since this is the only way a shorter string can be compatible with a longer one). And this means that $q$ must assign a new label to some node which lies along one of $y_1,\ldots,y_n$.

Now that we have established this property of $Y$, we can finish quickly. Define $p' = (g_p, X_p \cup Y)$. Note that by the definition of Kumabe-Slaman forcing, if a condition $q$ extends $p'$ then $q$ cannot assign any new labels along paths in $Y$. But this means that no such $q$ can be in $D$, contradicting the density of $D$.
\end{proof}

We will now explain how the key lemma above can be used to prove the property that we referred to earlier as the ``key property of Kumabe-Slaman forcing.'' Recall that this property was that if $M$ is a countable $\omega$-model of $\ZF$ and $x$ is a real not in $M$ then we can construct $g$ Kumabe-Slaman generic over $M$ such that $g(x)$ encodes essentially any information, including information about $g$ itself.

We will begin by explaining how to construct $g$ when the information to be encoded does not depend on $g$. More formally, we will prove the following proposition.

\begin{proposition}
\label{prop:encode1}
Suppose that $M$ is a countable $\omega$-model of $\ZF$, $x$ is any real not in $M$ and $a \in \Cantor$ is arbitrary. Then there is $g$ Kumabe-Slaman generic over $M$ such that $g(x) = a$.
\end{proposition}

\begin{proof}
Let $D_0, D_1, \ldots$ enumerate the dense sets for Kumabe-Slaman forcing over $M$. We will form a sequence of conditions $p_0 \geq p_1 \geq \ldots$ in $M$, alternately meeting dense sets in $M$ and coding more bits of $a$ into $g(x)$.

We will begin by setting $p_0 = (\0, \0)$---i.e.\ $p_0$ assigns no labels and has no forbidden paths. For each $n \geq 0$, first use Lemma~\ref{lemma:key} to extend $p_n$ to a condition $p_n'$ which meets $D_n$ but does not assign any new labels along $x$. Next, extend $p_n'$ to a condition $p_{n + 1}$ which codes the next bit of $a$. More formally, note that $x$ cannot be among the forbidden paths in $p_n'$ because $x$ is not even in $M$. Thus we can find an initial segment $\sigma$ of $x$ which is long enough that it does not lie along any of the finitely-many forbidden paths of $p_n'$ and which does not yet have any label assigned to it. Define $p_{n + 1}$ to be the condition extending $p_n'$ which assigns to $\sigma$ a label equal to the $n^\text{th}$ bit of $a$ and which explicitly decides not to label all other as-yet-unlabeled strings of length at most $|\sigma|$.

Now define $G$ to be the upwards closure of the $p_n$'s, i.e.\ $G = \{q \mid \exists n\, (p_n \leq q)\}$. It is straightforward to check that $G$ is a generic filter for Kumabe-Slaman forcing over $M$ and that if $g$ is the partial labelling associated to $G$ then $g(x) = a$.
\end{proof}

Using this same idea, it is possible to encode much more information into $g(x)$, including information that depends on $g$. In particular, if $G$ is the generic filter giving rise to $g$ then we may encode into $g(x)$ any fact about $M[G]$ which can be expressed in the forcing language over $M$. A precise statement is given in the proposition below. Note that by a \term{formula $\phi$ in the forcing language over $M$}, we mean a formula of the form $\phi(\dot{\gamma}_1, \ldots, \dot{\gamma}_n)$ where $\dot{\gamma}_1,\ldots,\dot{\gamma}_n$ are names in $M$ for Kumabe-Slaman forcing over $M$.

\begin{proposition}
\label{prop:encode2}
Suppose that $M$ is a countable $\omega$-model of $\ZF$, $x$ is any real not in $M$ and $\phi_0, \phi_1, \ldots$ is a sequence of formulas in the forcing language over $M$. Then there is $g$ Kumabe-Slaman generic over $M$ such that $g(x)$ encodes the truth of $\phi_0, \phi_1, \ldots$ in $M[G]$---i.e.\ the $n^\text{th}$ bit of $g(x)$ is $1$ if and only if $M[G] \models \phi_n$.
\end{proposition}

\begin{proof}
The proof is very similar to the proof of Proposition~\ref{prop:encode1}. That is, we let $D_0, D_1, \ldots$ enumerate the dense sets in $M$ and construct a sequence of conditions $p_0 \geq p_1 \geq \ldots$ by alternately meeting dense sets and encoding bits into $g(x)$. The only difference is that now, before encoding a bit of information into $g(x)$, we first have to determine whether the bit to be encoded should be a $0$ or $1$, which we can do by meeting an appropriate dense set.

Let's now be a bit more specific. For each formula $\phi_n$, fix a dense set $E_n$ of conditions which decide $\phi_n$---i.e.\ for each $p \in E_n$ we have either $p \forces \phi_n$ or $p \forces \lnot\phi_n$. Now set $p_0 = (\0, \0)$ and for each $n \geq 0$, extend $p_n$ to $p_{n + 1}$ as follows.

First, use Lemma~\ref{lemma:key} to extend $p_n$ to $p_n' \in D_n$ without assigning new labels along $x$. Next, use Lemma~\ref{lemma:key} to extend $p_n'$ to $p_n'' \in E_n$ without assigning new labels along $x$. Now we have that either $p_n'' \forces \phi$ or $p_n'' \forces \lnot \phi$. In the former case, extend $p_n''$ to $p_{n + 1}$ by assigning a new label of $1$ somewhere along $x$ (and not assigning any other new labels) and in the latter case, assign a label of $0$ somewhere along $x$.

We finish as in the proof of Proposition~\ref{prop:encode1} by setting $G = \{q \mid \exists n\, (p_n \leq q)\}$ and taking $g$ to be the associated partial labelling.
\end{proof}

Of course, we can also combine Propositions~\ref{prop:encode1} and~\ref{prop:encode2} in the obvious way: given any $a \in \Cantor$ and any sequence of formulas $\phi_0, \phi_1,\ldots$ in the forcing language over $M$, we can find $g$ Kumabe-Slaman generic over $M$ such that $g(x)$ computes both $a$ and the truth of each $\phi_n$ in $M[G]$.

\subsection{Kumabe-Slaman forcing preserves admissibility}

We will now use Lemma~\ref{lemma:key} from the previous section to prove that Kumabe-Slaman forcing preserves admissibility of countable ordinals. More precisely, we will prove the following theorem.

\begin{theorem}[Lutz-Siskind]
\label{thm:admissible}
Suppose $M$ is an $\omega$-model of $\ZF$ and $g$ is Kumabe-Slaman generic over $M$. Then for every real $x \in M$, $\omega_1^{x\oplus g} = \omega_1^x$.
\end{theorem}

This theorem was discovered in unpublished work by the author and Benjamin Siskind. The proof was inspired by ideas from Jananthan and Simpson~\cite{jananthan2021turing} (who were themselves fleshing out ideas sketched by Slaman).

In order to prove the theorem, we will make use of ill-founded $\omega$-models of $\ZF$. One piece of notation will be useful in discussing such models. Given an ill-founded $\omega$-model $M$, we will use $o(M)$ to denote the height of the well-founded part of $M$. Here's what we mean. If we consider the class of ordinals as interpreted in $M$, we can see that they form a linear order and that the ordinals of $M$ which are well-founded (i.e.\ which are isomorphic to ordinals in the real world) form a proper initial segment of this order. We will use $o(M)$ to refer to the height of this initial segment, i.e.
\[
  o(M) = \sup\{\alpha \in \Ord \mid \exists a \in \Ord^M \, (\text{$a$ is isomorphic to $\alpha$})\}.
\]
Note that $o(M)$ is itself an ordinal---in fact, it is always the least ordinal not isomorphic to any ordinal in $M$.

The first step in showing that Kumabe-Slaman forcing preserves admissibility is to show that admissibility is preserved when forcing over ill-founded $\omega$-models with sufficiently short well-founded parts.

\begin{lemma}
\label{lemma:preserves}
Suppose $x \in \Cantor$ is any real and $M$ is an ill-founded $\omega$-model of $\ZF$ such that $x \in M$ and $o(M) = \omega_1^x$. Then for every $g$ Kumabe-Slaman generic over $M$, $\omega_1^{x\oplus g} = \omega_1^x$.
\end{lemma}

\begin{proof}
Obviously $\omega_1^{x\oplus g} \geq \omega_1^x$; we need to show that $\omega_1^{g\oplus x} \leq \omega_1^x$. To that end, fix an index $e$ such that $\Phi_e(x\oplus g)$ is a presentation of an ordinal. We need to show that this ordinal is less than $\omega_1^x$.

Recall that $g$ is not literally a generic filter for Kumabe-Slaman forcing, but rather the partial labelling associated to one such filter. Let $G$ be the generic filter giving rise to $g$. Note that $g$ is definable from $G$ and hence $g$ is in $M[G]$. Thus we may talk about $\Phi_e(x\oplus g)$ in $M[G]$. Furthermore, since $M$ is an $\omega$-model, so is $M[G]$ and thus $M[G]$ is correct about values of $\Phi_e(x\oplus g)$. In other words, $\Phi_e(x\oplus g)$ means the same thing inside $M[G]$ that in means externally.

Since $\Phi_e(x\oplus g)$ is well-founded and $M[G]$ is an $\omega$-model, $M[G]$ believes $\Phi_e(x\oplus g)$ is well-founded. Thus it is isomorphic to some ordinal $\alpha$ in $M[G]$. And since $\Phi_e(x\oplus g)$ is really well-founded, $\alpha$ must be in the well-founded part of $M[G]$. But the ordinals of $M[G]$ are the same as those of $M$ and hence $\alpha < o(M) = \omega_1^x$.
\end{proof}

In Lemma~\ref{lemma:preserves} above, we essentially showed that Kumabe-Slaman forcing preserves admissibility of an ordinal $\alpha$ when we force over an ill-founded model $M$ with $o(M) = \alpha$. The next lemma will be used to show that to preserve admissibility of $\alpha$, it is enough to force over a model which \emph{contains} such an ill-founded model.

\begin{lemma}
\label{lemma:generic}
Suppose that $M$ is an $\omega$-model of $\ZF$ and $N \in M$ is an element of $M$ which $M$ believes to be an $\omega$-model of $\ZF$. Then for every $g$ which is Kumabe-Slaman generic over $M$, $g$ is also Kumabe-Slaman generic over $N$.
\end{lemma}

\begin{proof}
Note that since $M$ is an $\omega$-model and $M$ believes $N$ is an $\omega$-model, $N$ really is an $\omega$-model and thus we can apply Lemma~\ref{lemma:key} to it.

Let us now say more precisely what we need to prove. Let $G$ be a generic filter for Kumabe-Slaman forcing over $M$ for which $g$ is the associated partial labelling. We need to show that there is a generic filter $G'$ for Kumabe-Slaman forcing over $N$ such that $g$ is also the partial labelling associated to $G'$. It is perhaps tempting to think we can just take $G' = G$. The reason we can't (and the reason this lemma is not trivial) is that since possibly not all reals of $M$ are in $N$, possibly not all Kumabe-Slaman forcing conditions in $M$ are in $N$. Thus $G$ is not even a subset of $N$'s version of the Kumabe-Slaman forcing poset.

However, the next best thing does work: we can take $G'$ to be the intersection of $G$ with $N$. To see why this makes sense, note that since $N$ is an $\omega$-submodel of $M$, every Kumabe-Slaman condition in $N$ is also a condition in $M$. So we are justified in defining $G' = \{p \in G \mid p \in N\}$.

We first claim that $g$ is the partial labelling associated to $G'$. The reason is that, given any condition $p = (g_p, X_p)$ in $G$, the condition $p' = (g_p, \0)$ is also in $G$ and hence in $G'$.

We now claim that $G'$ is a generic filter over $N$. To prove this, we need to show that $G'$ is upwards closed, downwards directed and meets every dense set in $N$. The fact that $G'$ is upwards closed follows from the fact that $G$ is. To see that $G'$ is downwards directed, suppose that $p = (g_p, X_p)$ and $q = (g_q, X_q)$ are both in $G'$. Thus they are both in $G$, from which it is easy to see that so is $r = (g_p \cup g_q, X_p \cup X_q)$. And since $r$ is clearly in $N$, it is in $G'$ as well.

It remains to show that $G'$ meets every dense set in $N$. Let $D$ be such a dense set. We claim that $D$ generates a dense subset of the Kumabe-Slaman forcing poset in $M$. To be precise, define $D' = \{p \in M \mid \exists q \in D\, (p \leq q)\}$; we claim that $D'$ is dense in $M$.

To show this, fix $p = (g_p, X_p)$ a Kumabe-Slaman forcing condition in $M$. Let $p' = (g_p, \0)$. Since $N$ is an $\omega$-submodel of $M$, $p' \in N$. By Lemma~\ref{lemma:key}, we can find some $q = (g_q, X_q)$ extending $p'$ and in $D$ which does not assign any new labels along paths in $X_p$. This implies that $q' = (g_q, X_q \cup X_p)$ is a condition in $M$ extending both $p$ and $q$. Since $q'$ extends $q \in D$, $q' \in D'$, which shows that $D'$ is dense.

Since we know $D'$ is dense in $M$ and $G$ is generic, there must be some condition $p \in G$ which is in $D'$ and thus extends some condition $q \in D$. Since $p \leq q$ and $G$ is a filter, $q$ is also in $G$. Since $q$ is in $N$, $q$ is also in $G'$ and hence $G'$ meets $D$. Thus $G'$ is generic.
\end{proof}

We are now almost done with the proof of Theorem~\ref{thm:admissible}. To finish, it is essentially enough to show that every $\omega$-model contains enough ill-founded models. The following lemma will be used to prove this.

\begin{lemma}
\label{lemma:model}
For every real $x$, there is an ill-founded $\omega$-model $M$ of $\ZF$ such that $x \in M$ and $o(M) = \omega_1^x$.
\end{lemma}

\begin{proof}
Let $A_x$ be the set of reals which code countable $\omega$-models of $\ZF$ which contain $x$. Note that $A_x$ is $\Sigma^1_1(x)$ and nonempty. Thus by the Gandy basis theorem (see~\cite[Corollary III.1.5]{sacks1990higher}), there is an element $M$ of $A_x$ which is hyperlow relative to $x$---i.e.\ $\omega_1^{M \oplus x} = \omega_1^x$. We claim that $M$ has the desired properties.

First, by construction $M$ is an $\omega$-model of $\ZF$ which contains $x$. Second, since $M$ contains $x$ (and since $M$ is an $\omega$-model), $M$ must contain a copy of each ordinal computable from $x$. And since $M$ is an $\omega$-model, it must believe that these copies are all well-founded. Hence $o(M) \geq \omega_1^x$. It remains to show that $o(M) \leq \omega_1^x$. If we can show this then it follows that $M$ is ill-founded: $M$ believes that there are ordinals larger than $\omega_1^x$ and if $o(M) \leq \omega_1^x$ then these ordinals cannot actually be well-founded.

To show that $o(M) \leq \omega_1^x$, note that any ordinal in $M$ is computable from $M$. Thus if $\alpha$ is an ordinal in the well-founded part of $M$ then $\omega_1^M > \alpha$. By our choice of $M$, we have $\omega_1^M \leq \omega_1^{M\oplus x} = \omega_1^x$. So no ordinal in the well-founded part of $M$ can be greater than or equal to $\omega_1^x$.
\end{proof}

We can now prove Theorem~\ref{thm:admissible}.

\begin{proof}[Proof of Theorem~\ref{thm:admissible}]
This almost follows immediately from the preceding three lemmas, but there is a small problem. We will first explain the obvious, but slightly wrong, proof and then explain how to fix it.

Fix a real $x \in M$. Since $M$ is an $\omega$-model, $M$ knows that $x$ is a real. Since $M \models \ZF$, Lemma~\ref{lemma:model} holds in $M$ and thus $M$ contains some $N$ which $M$ believes to be an ill-founded $\omega$-model of $\ZF$ which contains $x$ and satisfies $o(N) = \omega_1^x$. It is not hard to check that this implies that $N$ really is an ill-founded $\omega$-model of $\ZF$ which contains $x$ and satisfies $o(N) = \omega_1^x$. Now let $g$ be Kumabe-Slaman generic over $M$. By Lemma~\ref{lemma:generic}, $g$ is also Kumabe-Slaman generic over $N$, so by Lemma~\ref{lemma:preserves}, $\omega_1^{x\oplus g} = \omega_1^x$.

The only problem with this proof is the claim that Lemma~\ref{lemma:model} holds in $M$. The problem is that this lemma cannot be proved in $\ZF$---it requires an assumption somewhat stronger than consistency of $\ZF$ to prove\footnote{The specific place in the proof where such an assumption is needed is for the assertion that $A_x$ is nonempty.}---and thus $M$ does not necessarily satisfy it despite being a model of $\ZF$. The solution is to observe that in the statements of all Lemmas, we can replace ``model of $\ZF$'' with ``model of some sufficiently large finite fragment of $\ZF$'' and that Lemma~\ref{lemma:model} is provable in $\ZF$ if we don't ask for a model of $\ZF$, but just for a model of any fixed finite fragment of $\ZF$.
\end{proof}

\bibliographystyle{alpha}
\bibliography{borel_graphable}

\end{document}